\documentclass[11pt]{amsart} \usepackage{amscd}
\usepackage{amsfonts,amssymb,latexsym} 
\input xy
\xyoption{all}
\xyoption{arc}
\xyoption{ps}
\xyoption{dvips}
\CompileMatrices

\newcommand{\tensor}{tensor}
\newcommand{\tensors}{tensors}
\newcommand{\Tensor}{Tensor}
\newcommand{\Tensors}{Tensors}
\newcommand{\fixed}{D}
\newcommand{\sfixed}{\mathcal{D}}
\newcommand{\param}{R}
\newcommand{\di}{n}

\newcommand{\SA}{{\mathcal{A}}}

\newcommand{\SC}{{\mathcal{C}}}
\newcommand{\SD}{{\mathcal{D}}}

\newcommand{\SF}{{\mathcal{F}}}
\newcommand{\SG}{{\mathcal{G}}}
\newcommand{\SH}{{\mathcal{H}}}
\newcommand{\SI}{{\mathcal{I}}}

\newcommand{\SK}{{\mathcal{K}}}

\newcommand{\SM}{{\mathcal{M}}}

\newcommand{\SO}{{\mathcal{O}}}
\newcommand{\SP}{{\mathcal{P}}}

\newcommand{\cS}{{\mathcal{S}}}
\newcommand{\ST}{{\mathcal{T}}}

\newcommand{\SV}{{\mathcal{V}}}

\newcommand{\PP}{\mathbb{P}}
\newcommand{\ZZ}{\mathbb{Z}}
\newcommand{\CC}{\mathbb{C}}
\newcommand{\QQ}{\mathbb{Q}}

\newcommand{\FM}{\mathfrak{M}}

\newcommand{\isom}{\cong}
\newcommand{\mi}{{\rm min}}

\newcommand{\Spec}{\operatorname{Spec}}

\newcommand{\codim}{\operatorname{codim}}

\newcommand{\Hom}{\operatorname{Hom}}

\newcommand{\Pic}{\operatorname{Pic}}
\newcommand{\Sym}{\operatorname{Sym}}
\newcommand{\id}{\operatorname{id}}

\newcommand{\surj}{\twoheadrightarrow}
\newcommand{\inj}{\hookrightarrow}
\newcommand{\too}{\longrightarrow}
\newcommand{\rk}{\operatorname{rk}}
\newcommand{\wt}{\widetilde}
\newcommand{\glv}{\operatorname{GL}(V)}
\newcommand{\slv}{\operatorname{SL}(V)}
\newcommand{\slw}{\operatorname{SL}(W)}
\newcommand{\glr}{\operatorname{GL}(r)}
\newcommand{\gln}{\operatorname{GL}(n)}
\newcommand{\orth}{\operatorname{O}(r)}
\newcommand{\sor}{{\operatorname{SO}(r)}}
\newcommand{\sympl}{\operatorname{Sp}(r)}
\newcommand{\glrc}{\operatorname{GL}(r,\CC)}
\newcommand{\orthc}{\operatorname{O}(r,\CC)}
\newcommand{\sorc}{{\operatorname{SO}(r,\CC)}}
\newcommand{\symplc}{\operatorname{Sp}(r,\CC)}

\newcommand{\ev}{{\rm ev}}

\newtheorem{proposition}{Proposition}[section]
\newtheorem{theorem}[proposition]{Theorem}
\newtheorem{definition}[proposition]{Definition}
\newtheorem{lemma}[proposition]{Lemma}

\newtheorem{corollary}[proposition]{Corollary}
\newtheorem{remark}[proposition]{Remark}

\numberwithin{equation}{section}

\title[Tensors and orthogonal sheaves]
{Stable tensors and moduli space of orthogonal sheaves}
\author[T. G\'omez, I. Sols]{Tom\'as L. G\'omez and Ignacio Sols}
\date{21 January 2003}

\address{T. G\'omez, School of Mathematics, Tata Institute of
Fundamental Research, Mumbai 400 005 (India)
(current address: Universidad Complutense de Madrid)}

\email{tomas@math.tifr.res.in, tgomez@alg.mat.ucm.es}

\address{I. Sols,
Departamento de Algebra, Facultad de Ciencias Matem\'aticas,
Universidad Complutense de Madrid, 28040 Madrid (Spain)}

\email{sols@mat.ucm.es}

\begin{document}

\begin{abstract}
Let $X$ be a smooth projective variety over $\CC$. 
We find the natural notion of semistable 
orthogonal bundle and construct the moduli space, which we compactify
by considering also orthogonal sheaves, i.e. pairs $(E,\varphi)$, 
where $E$
is a torsion free sheaf on $X$ and $\varphi$ is a symmetric 
nondegenerate (in the open set where $E$ is locally free) 
bilinear form on $E$. 
We also consider special orthogonal sheaves, by adding a
trivialization $\psi$ of the determinant of $E$ such
that $\det(\varphi)=\psi^2$;
and symplectic sheaves, by considering a 
form which is skewsymmetric.
More generally, we consider
semistable tensors, i.e. multilinear forms on a torsion free
sheaf, and construct their projective moduli space using GIT.
\end{abstract}

\maketitle

Let $X$ be a smooth projective variety of dimension $\di$ over $\CC$.
If $X$ is a curve, the moduli space of vector bundles 
was constructed by Mumford, 
Narasimhan and Seshadri.
If $\dim(X)>1$, to obtain a projective moduli space 
we have to consider also torsion free sheaves, and this was 
done by Gieseker, Maruyama and Simpson.

An \textit{orthogonal bundle} is a pair $(E,\varphi)$, where $E$ is a vector
bundle and
$$
\varphi:E\otimes E\too \SO_X
$$ 
is a bilinear symmetric nowhere degenerate form.
The nondegeneracy means that it induces an isomorphism
$E\to E^\vee$, hence $(\det E)^2\isom \SO_X$. 

A \textit{special orthogonal bundle} is a triple $(E,\varphi,\psi)$
where $E$ and $\varphi$ are as before, and 
$$
\psi:\det E \too \SO_X
$$
is an isomorphism such that $\det(\varphi)=\psi^2$ (this equation
means that for all points $x\in X$, if we choose a basis for
the fiber $E_x$, the determinant of the matrix associated to 
$\varphi$ at $x$ is equal to the square of the scalar associated
to $\psi$ at $x$).

A \textit{symplectic bundle} is a pair $(E,\varphi)$, where $E$ is a vector
bundle and $\varphi$ is a bilinear skewsymmetric nowhere degenerate
form.

Note that giving an orthogonal (or special orthogonal, or symplectic)
bundle is equivalent to giving a principal bundle with group
structure group $\orthc$ (or $\sorc$, or $\symplc$).
To obtain a projective moduli space we have to consider
also orthogonal, or special orthogonal or symplectic sheaves, i.e.
allowing $E$ to be a torsion free sheaf,
and then requiring $\varphi$ to be 
nondegenerate only on the open subset of $X$ where $E$ is
locally free. 
We say that a subsheaf $F$ of $E$ is isotropic
if $\varphi|_{F\otimes F}=0$.

An orthogonal, or special orthogonal or symplectic sheaf 
is called stable (respectively
semistable) if for all proper isotropic subsheaves $F$ of $E$
$$
P^{}_F + P_{F^\perp} \prec P^{}_E \quad
(\text{respectively}\;\preceq),
$$
where $P_F$ is the Hilbert polynomial of $F$, $F^\perp$ is the
sheaf perpendicular to $F$ with respect to the form $\varphi$,
and, as usual, the inequality between polynomials $P_1\prec P_2$
(respectively $\preceq$) means that $P_1(m)<P_2(m)$ (respectively
$\leq$) for 
$m\gg0$
(see sections \ref{secorthsheaves}
and \ref{secspecialgbundles} for 
precise definitions).

A similar problem was considered by Sorger \cite{So}. 
He works on a curve $C$ (not necessarily smooth) on a smooth
surface $S$, and  constructs
the moduli space of torsion free sheaves on $C$ together with
a symmetric form taking values on the dualizing sheaf $\omega_C$.
Faltings \cite{Fa} has considered principal bundles
on semistable curves. For $G$ orthogonal or symplectic he
considers a torsion free sheaf with a quadratic form, and
he also defines a notion of stability. For general
reductive group $G$ he uses the approach of loop groups.

More generally, we can consider triples $(E,\varphi,u)$ where
$E$ is a torsion free sheaf on $X$ and $\varphi$ is a 
non-zero homomorphism
$$
\varphi: (E^{\otimes s})^{\oplus c}   
\, \too \, (\det E)^{\otimes b} \otimes \fixed_u,
$$
where $\fixed_u$ is a locally free sheaf belonging to a fixed family
$\{\fixed_u\}_{u\in \param}$
parametrized by a scheme $\param$ (for instance, $\param$ could be
$\Pic^{a}(X)$, and then $\fixed_u$ is any line bundle with fixed
degree $a$,
or we can take $\param$ to be a point, and then $\fixed_u$ is a fixed
locally free sheaf). 
We call these
triples \textit{\tensors}. See section \ref{secstability} for the precise
definition. \Tensors\ generalize several objects that have
already appeared in the literature.
If $s=1$, $b=0$, $c=1$, and $\param$ is a point, these are the 
framed modules of Huybrechts and Lehn. They gave two constructions
of their moduli space: In \cite{H-L1} for $\dim(X)\leq 2$, based
in the ideas of Gieseker \cite{Gi}, and in \cite{H-L2} for arbitrary
dimension, following the ideas of Simpson \cite{Si}.
If $\dim(X)=1$, $s=2$, $b=0$, $c=1$, $\param$ is 
a point and $\fixed_u$ is a line bundle, these are the
conic bundles of \cite{G-S}.
If $\dim(X)=1$ and $\SD$ is a family of line bundles,  
these are the decorated
vector bundles whose moduli space was constructed by Schmitt
\cite{Sch}.

Using geometric invariant theory (GIT) as in \cite{Si}
and \cite{H-L2}, 
we construct the moduli space of
semistable \tensors\ (sections
\ref{secstability} to \ref{sectheorem}). 
This is used in sections
\ref{secorthsheaves} and \ref{secspecialgbundles}
to construct  the projective moduli space  
of classical sheaves.

Finally, in section \ref{secrelatedmodulis}, as a 
further application we obtain moduli spaces for 
\textit{$\glrc$-representational pairs}, 
i.e. pairs $(P,\sigma)$ consisting of a principal $\glrc$-bundle 
$P$ and a section $\sigma$ of the vector bundle associated to a 
fixed representation of $\glrc$. We can also 
consider a quasi-projective scheme $Y$ with an action of $\glrc$, 
and then we can take $\sigma$ 
to be a section of the associated fiber bundle with
fiber $Y$. Banfield \cite{Ba} and Mundet \cite{MR}
have given Hitchin-Kobayashi correspondences for these objects.

In a future paper we will consider principal $G$-bundles for 
more general groups.

\medskip
\noindent\textbf{Notation.}
If $f:Y \to Y'$ is a morphism, 
we denote $\overline{f} =\id_X\times f$.
If $E_S$ is a coherent sheaf on $X\times S$, we denote
$E_S(m):= E_S \otimes p^*_X\SO_X(m)$.
To simplify the notation, we will denote the complex groups 
$\glrc$, $\orthc$, $\symplc$,... by $\glr$, $\orth$, $\sympl$.
If $X$, $Y$, $Z$ are schemes, then $\pi^{}_X$, $\pi^{}_{Y\times Z}$, etc...
denote the corresponding projections from $X\times Y \times Z$.

If $P_1$ and $P_2$ are two polynomials, we write $P_1\prec P_2$ 
if $P_1(m)<P_2(m)$ for $m\gg 0$, and analogously for ``$\preceq$'' and
``$\leq$''.
We use the convention that whenever ``(semi)stable'' and
``$(\leq)$''  appear in a sentence, two statements
should be read: one with ``semistable'' and ``$\leq$'' and another 
with ``stable'' and ``$<$''.

\medskip
\noindent\textbf{Acknowledgments.}
We would like to thank U. Bhosle, N. Fakhruddin, M.S. Narasimhan, S. Ramanan, 
C.S. Seshadri and Ch. Sorger for discussions 
on this subject. 
The authors are members of VBAC (Vector Bundles on Algebraic Curves),
which is partially supported by EAGER (EC FP5 Contract no.
HPRN-CT-2000-00099) and by EDGE (EC FP5 Contract
no. HPRN-CT-2000-00101).
T.G. was supported by a postdoctoral fellowship of
Ministerio de Educaci\'on y Cultura (Spain).

\section{Stability of \tensors}
\label{secstability}

Let $X$ be a smooth projective 
variety over $\CC$. Fix an ample line bundle
$\SO_X(1)$ on $X$. 
Fix a polynomial
$P$ of degree $n=\dim(X)$, and integers $s$, $b$, $c$. 
We will denote by $r$ and $d$ the rank and degree of a sheaf with
Hilbert polynomial $P$.
Fix a family $\{\fixed_u\}_{u\in \param}$ of locally free sheaves $X$
parametrized
by a scheme $\param$, i.e. we fix a locally free sheaf $\sfixed$ on 
$X\times \param$, and given a point $u\in \param$, 
we denote by $\fixed_u$ the restriction
to the slice $X\times u$.

\begin{definition}[\Tensor]
A \tensor\ is a triple $(E,\varphi,u)$, where $E$ is a coherent sheaf on $X$
with Hilbert polynomial $P_E=P$, $u$ is a point in $\param$, 
and $\varphi$ is a homomorphism
$$
\varphi: (E^{\otimes s})^{\oplus c}   
\, \too \, (\det E)^{\otimes b} \otimes \fixed_u,
$$
that is not identically zero.
Let $(E,\varphi,u)$ and $(F,\psi,v)$ be two \tensors\ with $P_E=P_F$,
$\det E\cong \det F$, and $u=v$.
A homomorphism between $(E,\varphi,u)$ and $(F,\psi,v)$
is a pair $(f,\alpha)$ where $f:E\too F$ is a homomorphism of sheaves,
$\alpha\in \CC$,  and 
the following diagram commutes
\begin{equation}
\label{isomorphism}
\xymatrix{
{(E^{\otimes s})^{\oplus c}} 
\ar[rr]^{(f^{\otimes s})^{\oplus c}} \ar[d]^{\varphi} & & 
{(F^{\otimes s})^{\oplus c}} \ar[d]^{\psi}\\
{(\det E)^{\otimes b}\otimes\fixed_u} 
\ar[rr]^{\hat{f}\otimes \alpha} & & 
{(\det F)^{\otimes b}\otimes\fixed_v}
}
\end{equation}
where $\hat{f}:\det E\to \det F$ is the homomorphism induced by $f$.
In particular, $(E,\varphi,u)$ and $(E,\lambda\varphi,u)$ are isomorphic
for $\lambda\in \CC^*$.
\end{definition}

\begin{remark}
\label{alpha1}
\textup{
We could have defined a more restrictive notion of isomorphism,
considering only isomorphisms for which $\alpha=1$. If we do 
this,  we obtain a different category: for instance, if 
$E$ is simple, the set of automorphisms of $(E,\varphi,u)$ is
$\CC^*$, but if we require $\alpha=1$, then the set of
automorphisms is $\ZZ/(rb-s)\ZZ$ (assuming $rb-s\neq 0$).
If $rb-s\neq 0$, even if the categories are not equivalent, 
the set of isomorphism classes will be the same (because
$\alpha$ can be absorbed in $f$ by changing $f$ into
$\alpha^{1/(rb-s)}f$), and then the moduli spaces will 
be the same.
But if $rb-s=0$, then $\alpha$ cannot be absorbed in $f$,
and the set of isomorphism classes is not the same.}
\end{remark}

Let $\delta$ be a polynomial with $\deg(\delta)<\di=\dim(X)$
\begin{equation}
\label{delta}
\delta = \delta_1 t^{\di-1} + \delta_2 t^{\di-2} +\dots + \delta_\di \;\in
\QQ[t],
\end{equation}
and $\delta(m)> 0$ for $m\gg 0$. We denote $\tau=\delta_1 (\di-1)!$.
We will define a notion of stability for
these \tensors, depending on the
polarization and $\delta$, and we will construct, using geometric
invariant theory (GIT), a moduli space
for semistable \tensors.

A weighted filtration $(E_\bullet,m_\bullet)$ of 
a sheaf $E$ is a filtration of sheaves
\begin{equation}
\label{filtE}
0 \subsetneq E_1 \subset E_2 \subset \;\cdots\; \subset E_t \subsetneq
E_{t+1}=E,
\end{equation}
and positive numbers $m_1,\, m_2,\ldots , \,m_t>0$.
Let $r_i=\rk(E_i)$. If $t=1$ (one step filtration), then we will take
$m_1=1$. We will denote $E^i=E/E_i$ and $r^i=\rk(E^i)$.
The filtration is called \textit{saturated} if all
sheaves $E_i$ are saturated in $E$, i.e. if $E^i$ is torsion free.

Consider the vector of $\CC^r$ defined as $\gamma=\sum_{i=1}^{t} 
m_i \gamma^{(r_i)}$, where
\begin{equation}
\label{gammak}
\gamma^{(k)}=\big( \overbrace{k-r,\dots,k-r}^k,
\overbrace{k,\dots,k}^{r-k} \big)
\qquad
(1\leq k \leq  r-1).
\end{equation}
Now let $\SI=\{1,\ldots,t+1\}^{\times s}$ be the set of all
multi-indexes $I=(i_1,\ldots,i_s)$. Let $\gamma_{j}$ be the 
$j$-th component of the vector $\gamma$, and define
\begin{equation}
\label{muE}
\mu(\varphi,E_\bullet,m_\bullet)=\min_{I\in \SI} \big\{
\gamma_{r_{i_1}}+\dots+\gamma_{r_{i_s}}: \,
\varphi|_{(E_{i_1}\otimes\cdots \otimes E_{i_s})^{\oplus c}}\neq 0
  \big\}
\end{equation}

\begin{definition}[Stability]
\label{stability}
Let $\delta$ be a polynomial as in (\ref{delta}).
We say that $(E,\varphi,u)$ is $\delta$-(semi)stable if for all weighted
filtrations it is
\begin{equation}
\label{stformula}
\Big(\sum_{i=1}^t m_i\big( r P_{E_i} -r_i P \big)\Big) +
\mu(\varphi,E_\bullet,m_\bullet) \, \delta \;(\preceq)\; 0
\end{equation}
\end{definition}

Recall that we assume that $\varphi$ is not identically zero. 
It is easy to check that if $(E,\varphi,u)$ is semistable, then
$E$ is torsion free. In this definition, it suffices to consider
saturated filtrations, and with $\rk(E_i)<\rk(E_{i+1})$ for all $i$.

\begin{lemma}
\label{finiteE}
There is an integer $A_1$ (depending only on $P$, $s$, $b$, $c$ and
$\sfixed$) such that it is enough to check the stability condition 
(\ref{stformula}) for
weighted filtrations with $m_i\leq A_1$ for all $i$.
\end{lemma}

\begin{proof}
Again, let $\SI=\{1,\ldots,t+1\}^{\times s}$. Multi-indexes will be
denoted $I=(i_1,\ldots,i_s)$.
Note that (\ref{muE}) is a piece-wise linear function of 
$\gamma\in \SC$, where $\SC\subset \ZZ^r$ is the cone defined
by $\gamma_1\leq\ldots\leq\gamma_r$. This is because it is defined as
the minimum among a finite set of linear functions, namely the functions
$\gamma_{r_{i_1}}+\cdots+\gamma_{r_{i_s}}$ for 
$I\in \SI$. 
There is a decomposition of $\SC=\bigcup_{I\in \SI} \SC_I$ 
into a finite number of subcones
$$
\SC_I := \big\{ \gamma\in \SC :
\gamma_{r_{i_1}}+\cdots+\gamma_{r_{i_s}} \leq
\gamma_{r_{i'_1}}+\cdots+\gamma_{r_{i'_s}} \;\text{for all } I'
\in \SI \big \} 
$$
and  (\ref{muE}) is linear on each cone $\SC_I$.
Choose one vector $\gamma\in
\ZZ^r$ in each edge of each cone $\SC_I$. 
Multiply all these vectors by $r$, so that all their
coordinates are divisible by $r$, and call this set of vectors $S$. 
All vectors in $S$ 
come from a collection of weights
$m_i>0$, $i=1,\dots,t+1$, given by the formula 
$\gamma=\sum_{i=1}^{t} m_i \gamma^{(r_i)}$. Hence to obtain the finite
set $S$ of vectors it is enough to consider a finite set of values for
$m_i$, and hence there is a maximum value $A_1$.

Finally, it is easy to see that it is enough to check
(\ref{stformula}) for the weights associated to the vectors in $S$. 
Indeed, 
first note that since the first term in (\ref{stformula}) is linear on
$\SC$, then it is also linear on each $\SC_I$. Then the left hand side
of (\ref{stformula}) is linear on each $\SC_I$, and hence it is enough
to check it on all the edges of all the cones $\SC_I$.

\end{proof}

\begin{definition}[Slope stability]
\label{slopestability}
We say that $(E,\varphi,u)$ is slope-$\tau$-(semi)stable if 
$E$ is torsion free, and for all weighted
filtrations we have
\begin{equation}
\Big(\sum_{i=1}^t m_i\big( r \deg E_i  -r_i \deg E \big)\Big) +
\mu(\varphi,E_\bullet,m_\bullet) \, \tau \;(\leq)\; 0
\end{equation}
\end{definition}

Recall. $\tau=\delta_1 (\di-1)!$.
As usual, we have the following implications
$$
\text{slope-$\tau$-stable $\Longrightarrow$ $\delta$-stable
$\Longrightarrow$}
$$
$$
\text{$\Longrightarrow$ $\delta$-semistable $\Longrightarrow$ 
slope-$\tau$-semistable}
$$

The reason why we have to consider filtrations 
instead of just
subsheaves is that (\ref{muE}) is not linear as a function of
$\{m_i\}$. But we have the following result that will be used in the
proof of theorem \ref{thm2.1}.

\begin{lemma}
\label{estimate}
Let $(E_\bullet,m_\bullet)$ be a filtration as above, and let $\ST'$ be a
subset of $\ST=\{1,\ldots,t\}$. Let
$(E'_\bullet,m'_\bullet)$ be the subfiltration
obtained by considering only those terms 
$E_i$ for which $i\in \ST'$.
Then
$$
\mu(\varphi,E_\bullet,m_\bullet) \;\leq\;
\mu(\varphi,E'_\bullet,m'_\bullet)
+\sum_{i\in \ST-\ST'} m_i \, s \, r_i \,.
$$
\end{lemma}

\begin{proof}
Let $\SI=\{1,\ldots,t+1\}^{\times s}$ be the set of all
multi-indexes $I=(i_1,\ldots,i_s)$.
Given a multi-index $I\in \SI$, we have
$$
\gamma_{r_{i_1}}+\dots+\gamma_{r_{i_s}}=\sum_{i=1}^t
m_i\big( s\,r_i-\nu_i(I)\,r \big)\, ,
$$
where $\nu_i(I)$ is the number of elements $k$ of the multi-index
$I=(i_1,\ldots,i_s)$ such that $r_k\leq r_i$. 
If $I$ is
the multi-index giving minimum in (\ref{muE}), we will denote
$\epsilon_i(\varphi,E_\bullet,m_\bullet):=\nu_i(I)$
(or just $\epsilon_i(E_\bullet)$ if the rest of the data is
clear from the context). Then
\begin{equation}
\label{epsilonE}
\mu(\varphi,E_\bullet,m_\bullet)=
\sum_{i=1}^t m_i \big( s\,r_i-\epsilon_i(E_\bullet)\,r \big)\, .
\end{equation}
We index the filtration $(E'_\bullet,m'_\bullet)$ with $\ST'$.
Let $I'=(i'_1,\dots,i'_s)\in \big\{\ST'\cup \{t+1\}\big\}^{\times s}$ 
be the multi-index giving minimum for the filtration
$(E'_\bullet,m'_\bullet)$. In particular, we have 
$\varphi|_{(E_{i'_1}\otimes\cdots \otimes E_{i'_s})^{\oplus c}}\neq 0$.
Then
\begin{eqnarray*}
&&\mu(\varphi,E_\bullet,m_\bullet)  \;=\; 
      \min_{I\in \SI} \big\{
      \gamma_{r_{i_1}}+\dots+\gamma_{r_{i_s}}: \,
      \varphi|_{(E_{i_1}\otimes\cdots \otimes E_{i_s})^{\oplus c}}\neq 0
      \big\}   \\
&&\hspace{1cm} \leq\;  \gamma_{r_{i'_1}}+\dots+\gamma_{r_{i'_s}}    \\
&&\hspace{1cm} =\; \sum_{i=1}^t m_i \big( s\,r_i-\nu_i(I')\,r \big) \\
&&\hspace{1cm} =\; \sum_{i=1}^t m_i 
       \big( s\,r_i-\epsilon_i(E'_\bullet)\,r \big) \\
&&\hspace{1cm} =\; \sum_{i\in \ST'} m_i 
      \big( s\,r_i-\epsilon_i(E'_\bullet)\,r \big)
      \;+\; \sum_{i\in \ST-\ST'} m_i 
      \big( s\,r_i-\epsilon_i(E'_\bullet)\,r \big) \\
&&\hspace{1cm} \leq\; \mu(\varphi,E'_\bullet,m'_\bullet)
      \;+\; \sum_{i\in \ST-\ST'} m_i\, s\, r_i\,. 
\end{eqnarray*}
\end{proof}

A family of $\delta$-(semi)stable \tensors\ parametrized by a scheme $T$ is a
tuple $(E_T,\varphi_T,u_T,N)$, 
consisting of a torsion free sheaf $E_T$ on 
$X\times T$, flat over $T$, that restricts to a torsion free sheaf
with Hilbert polynomial $P$ on every slice $X\times t$, 
a morphism $u_T:T\to \param$, a line bundle $N$ on $T$ and
a homomorphism $\varphi_T$  
\begin{equation}
\label{family}
\varphi^{}_T:(E^{}_T{}^{\otimes s})^{\oplus c}   
\, \too \, (\det E^{}_T)^{\otimes b} \otimes 
\overline{u^{}_T}^* \sfixed \otimes \pi_T^*N,
\end{equation}
(recall $\overline{u_T}=\id_X\times u_T$) 
such that if we consider the restriction of this
homomorphism on every slice $X\times t$ 
$$
\varphi^{}_t:(E^{}_t{}^{\otimes s})^{\oplus c}   
\, \too \, (\det E^{}_t)^{\otimes b} \otimes \fixed^{}_{u(t)},
$$
the triple $(E_t,\varphi_t,u(t))$
is a 
$\delta$-(semi)stable \tensor\ for every $t$ 
(in particular, $\varphi_t$ is not
identically zero).

Two families $(E_T,\varphi_T,u_T,N)$ and
$(E'_T,\varphi'_T,u_T',N')$ parametrized by $T$ 
are isomorphic if $u_T=u_T'$ and there are isomorphisms
$f:E_T\to E_T'$, 
$\alpha:N\to N'$ such that the induced diagram
$$
\xymatrix{
{(E_T^{}{}^{\otimes s})^{\oplus c}} 
\ar[rr]^{(f^{\otimes s})^{\oplus c}} \ar[d]^{\varphi_T} & & 
{({E_T'}{}^{\otimes s})^{\oplus c}} \ar[d]^{\varphi_T'}\\
{(\det E_T)^{\otimes b}\otimes \overline{u_T}^* \sfixed
\otimes \pi_T^*N} 
\ar[rr]^{\hat{f}\otimes\pi_T^*\alpha} & & 
{(\det E_T')^{\otimes b}\otimes \overline{u_T'}^* \sfixed
\otimes \pi_T^*N'}
}
$$
commutes.

Let ${\SM^{}_\delta}$ (respectively
${\SM^s_\delta}$) be the contravariant 
functor from the category of complex schemes, locally of
finite type, $({\rm Sch}/\CC)$ to the category of sets
$({\rm Sets})$ which sends a scheme
$T$ to the set of isomorphism classes of families of $\delta$-semistable
(respectively stable) \tensors\ parametrized by $T$, and sends a 
morphism $T'\to T$ to the map defined by pullback (as usual).

\begin{definition}
\label{corepresents}
A scheme $Y$ corepresents a functor $F:({\rm Sch}/\CC)
\to ({\rm Sets})$ if
\begin{enumerate}

\item There exists a natural transformation $f:F \to \underline{Y}$
(where $\underline{Y}$ is the functor of points represented by $Y$).

\item For every scheme $Y'$ and natural transformation 
$f':F \to \underline{Y'}$, there exists a unique $g:\underline{Y}\to
\underline{Y'}$ such that $f'$ factors through $f$.

\end{enumerate}
If $Y$ exists, then it is unique up to unique isomorphism. If 
furthermore $f(\Spec \CC):F(\Spec\CC) \to Y$ is bijective,
we say that $Y$ is a coarse moduli space. 
\end{definition}

We will construct schemes $\FM^{}_\delta$, $\FM^s_\delta$ corepresenting
the functors ${\SM^{}_\delta}$ and 
${\SM^s_\delta}$. In general 
$\FM_\delta$
will not be a coarse moduli space, because nonisomorphic \tensors\
could
correspond to the same point in $\FM_\delta$. 
As usual, we declare 
two such \tensors\ S-equivalent, and then $\FM_\delta$  becomes a coarse
moduli space for the functor of S-equivalence classes of \tensors.

\begin{theorem}
\label{mainthm}
Fix $P$, $s$, $b$, $c$ and a family $\SD$ of locally free sheaves
on $X$ parametrized by a scheme $\param$.
Let $d$ be the degree of a coherent sheaf whose Hilbert polynomial
is $P$. Let $\delta$ be a polynomial as in (\ref{delta}).
 
There exists a coarse moduli space $\FM_\delta$, projective over
$\Pic^d(X)\times \param$, of
S-equivalence classes of $\delta$-semistable 
\tensors.
The closed points of $\FM_\delta$ correspond to
S-equivalence classes of $\delta$-semistable \tensors. 
There is an open set $\FM^s_\delta$ 
corresponding to $\delta$-stable \tensors. 
Points in this open set
correspond to isomorphism classes of $\delta$-stable \tensors.
\end{theorem}

In proposition \ref{sequiv} we give a criterion to decide when 
two \tensors\ are S-equivalent.
Theorem \ref{mainthm} will be proved in section \ref{sectheorem}.

\begin{remark}\textup{ 
Note that to define the functors 
we have used \textit{isomorphism} classes of families,
but usually one uses \textit{equivalence} classes, declaring
two families equivalent if they differ by the pullback of a
line bundle $M$ on $T$. 
As a result, in general the functors that
we have defined will not be sheaves. The sheafified
functors will be the same (this follows from the fact
that if we shrink $T$ then $M$ will be trivial),
and hence the corresponding
moduli spaces are the same,
because a scheme corepresents a functor if and only if it
corepresents its sheafification (see \cite[p. 60]{Si}).
}
\end{remark}

\section{Boundedness}
\label{secboundedness}

The objective of this section is theorem \ref{thm2.1},
where we reformulate the stability condition for \tensors\ 
using some boundedness results.
We start with some well known results. See \cite[cor 1.7]{Si} (also
\cite[lemma 2.2]{H-L2}), \cite[lemma 2.5]{Gr} and \cite{Ma}.

\begin{lemma}[Simpson]
\label{lem2.2}
Let $r>0$ be an integer. Then there exist a constant $B$ with the
following property: for every torsion free 
sheaf $E$ with $0<\rk(E)\leq r$, we have
$$
h^0(E)\;\leq\;
\frac{1}{g^{\di-1} \di!}
\Big( \big( \rk(E)-1\big) \big([\mu_{\rm max}(E)+B]_+\big)^\di
+ \big([\mu_{\rm min}(E)+B]_+\big)^\di \Big),
$$
where $g=\deg \SO_X(1)$, $[x]_+=\max\{0,x\}$, and $\mu_{\rm max}(E)$
(respectively $\mu_{\rm min}(E)$) is the maximum (respectively minimum) slope of the
Mumford-semistable factors of the Harder-Narasimhan filtration of
$E$. 
\end{lemma}

\begin{lemma}[Grothendieck]
\label{lem2.6}
Let $\ST$ be a bounded set of sheaves $E$. The set of 
torsion free quotients $E''$ of the sheaves $E$ in $\ST$
with $|\deg(E'')|\leq C''$ for some fixed constant $C''$, 
is bounded.
\end{lemma}

\begin{theorem}[Maruyama]
\label{thm2.3}

The family of sheaves $E$ with fixed Hilbert polynomial $P$ and such
that $\mu_{\rm max}(E)\leq C$ for a fixed constant $C$, is bounded.
\end{theorem}

\begin{corollary}
The set of $\delta$-semistable \tensors\ $(E,\varphi,u)$ 
with fixed Hilbert polynomial is bounded.
\end{corollary}

\begin{proof}
Follows from theorem \ref{thm2.3} and an easy calculation.
\end{proof}

The main theorem of this section is

\begin{theorem}
\label{thm2.1}
There is an integer $N_0$ such that if $m\geq N_0$, the following
properties of \tensors\ $(E,\varphi,u)$ with $E$ torsion free 
and $P_E=P$, are equivalent.

\begin{enumerate}

\item $(E,\varphi,u)$ is (semi)stable.

\item $P(m)\leq h^0(E(m))$ and for every weighted filtration 
$(E_\bullet,m_\bullet)$ as in
(\ref{filtE})
$$
\Big(\sum_{i=1}^t m_i\big(r\, h^0(E_i(m))-r_i P(m)\big) \Big)+
\mu(\varphi,E_\bullet,m_\bullet)\, \delta(m) \;(\leq)\; 0
$$

\item For every weighted filtration $(E_\bullet,m_\bullet)$ 
as in (\ref{filtE})
$$
\Big(\sum_{i=1}^t m_i\big(r^i P(m)-r h^0(E^i(m))\big) \Big) +
\mu(\varphi,E_\bullet,m_\bullet)\, \delta(m) \;(\leq)\; 0
$$

\end{enumerate}

Furthermore, for any \tensor\ $(E,\varphi,u)$ satisfying these conditions, 
E is $m$-regular.

\end{theorem}

Recall that a sheaf $E$ is called $m$-regular if $h^i(E(m-i))=0$ for
$i>0$. If $E$ is $m$-regular, then $E(m)$ is generated by global sections,
and it is $m'$-regular for any $m'>m$.
The set of \tensors\ $(E,\varphi,u)$, with $E$ torsion free and $P_E=P$,
satisfying the weak version of conditions 1-3 will be called
$\cS^s$, $\cS'_m$, and $\cS''_m$.

\begin{lemma}
\label{lem2.4}
There are integers $N_1$, $C$ such that if $(E,\varphi,u)$ belongs to 
$\cS=\cS^s \cup \bigcup_{m\geq N_1} \cS''_m$, then for all 
saturated weighted filtrations the following holds for all $i$:
\begin{equation}
\label{2.4.5}
\deg(E_i)-r_i\mu_s \leq C,
\end{equation}
(where $\mu_s=(d-s\tau)/r$) 
and either $-C \leq \deg(E_i)-r_i \mu_s$ or

\begin{enumerate}

\item $r\, h^0(E_i(m))< r_i(P(m)-s\delta(m))$, if $(E,\varphi,u)\in \cS^s$
and $m\geq N_1$.

\item $r^i (P-s\,\delta) \prec r(P_{E^i}-s\,\delta)$, 
if $(E,\varphi,u)\in
\bigcup_{m\geq N_1} \cS''_m$.

\end{enumerate}
\end{lemma}

\begin{proof}
Let $B$ be as in lemma \ref{lem2.2}. Choose $C$ large enough so that 
$C>s\tau$ and the leading coefficient of $G-(P-s\delta)/r$ is
negative, where 
\begin{equation}
\label{2.4.2}
G(m)=\frac{1}{g^{\di-1}\di!}\Big( \big(1-\frac{1}{r}\big)
\big(\mu_s+s\tau+mg+B\big)^\di
+\frac{1}{r}\big(\mu_s-\frac{1}{r}C+mg+B\big)^\di \Big)
\end{equation}
Choose $N_1$ large enough so that for $m\geq N_1$, 
\begin{eqnarray}
\delta(m)&\geq&0, \label{2.4.3} \\ 
\mu_s-\frac{C}{r}+mg+B&>&0, \label{2.4.3b} \\
G(m)-(P(m)-s\delta(m))/r&<&0. \label{2.4.4}
\end{eqnarray}
Since the filtration is assumed to be saturated, and since $E$ is
torsion free, we have $0<r_i<r$.

\noindent\textbf{Case 1.} Suppose $(E,\varphi,u)\in \cS^s$. For each
$i$, consider the one step filtration $E_i \subsetneq E$. 
The leading coefficient of the semistability condition
applied to this filtration, together with $C>s\tau$, implies
(\ref{2.4.5}).

Let $E_{i,\max}\subset E_i$ be the term in the Harder-Narasimhan
filtration with maximal slope. Then the same argument applied to
$E_{i,\max}$ gives
\begin{equation}
\label{2.4.6}
\mu_{\max}(E_i)\,=\,\mu(E_{i,\max})\,<\, \mu_s + s\,\tau.
\end{equation}
Now assume that the first alternative does not hold, i.e.
$$
-C \,>\, \deg(E_i)-r_i \mu_s.
$$
This gives 
\begin{equation}
\label{2.4.7}
\mu_\mi(E_i) \,\leq\, \mu(E_i) \,<\, \mu_s-\frac{C}{r}.
\end{equation}
Combining lemma \ref{lem2.2} with (\ref{2.4.3b}), (\ref{2.4.6}), 
(\ref{2.4.7}) and 
(\ref{2.4.4}), we have 
$$r \, h^0(E_i(m)) \,<\, r_i(P(m)-s\,\delta(m)).$$

\noindent\textbf{Case 2.} Suppose $(E,\varphi,u)\in \cS''_m$ for some
$m\geq N_1$.
For each $i$, consider the quotient $E^i=E/E_i$. Let $E^i_\mi$ be the
last factor of the Harder-Narasimhan filtration of $E^i$ (i.e. 
$\mu(E^i_\mi)=\mu_\mi(E^i)$). Let $E'$ be the kernel
$$
0 \too E' \too E \too E^i_\mi \too 0,
$$
and consider the one step filtration $E'\subsetneq E$.
Equations (\ref{2.4.2}) and (\ref{2.4.3b}) imply that $0<G(m)$. Then
a short calculation using (\ref{2.4.4}), the fact that
$(E,\varphi,u)\in \cS''_m$, (\ref{2.4.3}) and lemma \ref{lem2.2} shows
$$
G(m) \;<\; \frac{h^0(E^i_\mi(m))}{\rk(E^i_\mi)} \;\leq\;
\frac{1}{g^{\di-1}\di!}\big(\mu_\mi(E^i)+mg+B\big)^\di.
$$
It can be seen that if this inequality of polynomials 
holds for some $m\geq N_1$, 
then it holds for all larger values of $m$, 
hence choosing $m$ large enough and looking at the coefficients, 
we have
$$
\mu_\mi(E^i) \,\geq\, \mu_s +\big(1-\frac{1}{r}\big)\,s\,\tau -\frac{C}{r^2}.
$$
A short calculation using this, 
$\mu_\mi(E^i)\leq\mu(E^i)$ and $0<\rk(E^i)<\rk(E)$ (hence
$\rk(E)>1$), yields (\ref{2.4.5}).

Now assume that the first alternative does not hold, i.e.
$$
-C \,>\, \deg(E_i)-r_i \mu_s.
$$
It follows that $r^i \mu_s < \deg(E^i)-s\,\tau$,
and hence
$$
r^i (P-s\delta) \,\prec\, r(P_{E^i}-s\delta).
$$
\end{proof}

\begin{lemma}
\label{lem2.5}
The set $\cS=\cS^s \cup \bigcup_{m\geq N_1} \cS''_m$ is bounded.
\end{lemma}

\begin{proof}
Let $(E,\varphi,u)\in \cS$. Let $E'$ be a subsheaf of $E$, and
$\overline{E}'$ the saturated subsheaf of $E$ generated by $E'$.
Using lemma \ref{lem2.4}
$$
\frac{\deg (E')}{\rk(E')} \;\leq\;
\frac{\deg (\overline{E}')}{\rk(\overline{E}')}
\;\leq\; \mu_s +\frac{C}{\rk(\overline{E}')} \;\leq\; \mu_s +C.
$$
Then by Maruyama's theorem \ref{thm2.3}, the set $\cS$ is bounded.
\end{proof}

\begin{lemma}
\label{lem2.7}
Let $\cS_0$ be the set of sheaves $E'$ such that $E'$ is a saturated
subsheaf of $E$ for some $(E,\varphi,u)\in \cS$, and furthermore
$$
|\deg(E')-r'\mu_s| \;\leq\; C.
$$
Then $\cS_0$ is bounded.
\end{lemma}

\begin{proof}
Let $E'\in \cS_0$. The sheaf $E''=E/E'$ is torsion free, and
$\;|\deg(E'')|\;$ is bounded because the set $\cS$ is bounded and
$$
|\deg(E'')| \;\leq\; |\deg(E)|+|\deg(E')| \;\leq\; \max_{E\in \cS}{|\deg (E)|}
+ C + r|\mu_s|.
$$
Then by Grothendieck's lemma \ref{lem2.6}, the set of sheaves $E''$
obtained in this way is bounded, and hence also $\cS_0$.
\end{proof}

\begin{lemma}
\label{lem2.8}
There is an integer $N_2$ such that for every weighted filtration 
$(E_\bullet,m_\bullet)$ as in (\ref{filtE}) 
with $E_i\in \cS_0$, the inequality
of polynomials
$$
\Big( \sum_{i=1}^t m_i \big( r P_{E_i}-r_i P \big) \Big)
+\mu(\varphi,E_\bullet,m_\bullet)\,\delta \;(\preceq)\; 0
$$
holds if and only if it holds for a particular value of $m\geq N_2$.
\end{lemma}

\begin{proof}
Since $\cS_0$ is bounded, the set that consists of the polynomials
$\delta$, $P_0$, $r'P_0$ and $P_{E'}$ for $E'\in \cS_0$ is finite. 
On the other hand, lemma \ref{finiteE} implies that we only need to
consider a finite number of values for $m_i$,
hence the result follows. 
\end{proof}

\begin{proof}[Proof of theorem \ref{thm2.1}]
Let $N_0>\max\{N_1,N_2\}$ and such that all sheaves in $\cS$ and
$\cS_0$ are $N_0$-regular, and $E_{1}\otimes \cdots \otimes E_{s}$
is $sN_0$-regular for all $E_{1}, \ldots , E_{s}$ in $\cS_0$.

\medskip
\noindent\textbf{2. $\Rightarrow$ 3.}
Let $(E,\varphi,u)\in \cS'_m$. Consider a weighted filtration
$(E_\bullet,m_\bullet)$. Then
\begin{eqnarray*}
\Big(\sum_{i=1}^t m_i(r^i P(m)-r h^0(E^i(m))) \Big) +
\mu(\varphi,E_\bullet,m_\bullet)\, \delta(m) &\leq \\
\Big(\sum_{i=1}^t (r h^0(E_i(m))-r_i P(m)) \Big)+
\mu(\varphi,E_\bullet,m_\bullet)\, \delta(m) &(\leq)& 0
\end{eqnarray*}

\medskip
\noindent\textbf{1. $\Rightarrow$ 2.}
Let $(E,\varphi,u)\in \cS^s$ and consider a saturated weighted filtration 
as in (\ref{filtE}). Since $E$ is $N_0$-regular, $P(m)=h^0(E(m))$.
If $E_i\in \cS_0$, then $P_{E_i}(m)=h^0(E_i(m))$.
If $E_i\notin \cS_0$, then the second alternative of lemma
\ref{lem2.4} holds, and then
$$
r \, h^0(E_i(m)) \;<\; r_i(P(m)-s\,\delta(m))\, .
$$
Let $\ST'\subset \ST=\{1,\ldots,t\}$ be the subset of those $i$
for which $E_i\in \cS_0$. Let $(E'_\bullet,m'_\bullet)$ be the
corresponding subfiltration.
Lemma \ref{estimate} and a short 
calculation shows that
\begin{eqnarray}
\label{chain}
\Big(\sum_{i=1}^t m_i \big(r \,h^0(E_i(m))-r_i P(m) \big) \Big) +
\mu(\varphi,E_\bullet,m_\bullet) \,\delta(m) & \leq& \\
\nonumber
\Big(\sum_{i\in \ST'} m_i \big(r \,P_{E_i}(m)-r_i P(m) \big) \Big) +
\mu(\varphi,E'_\bullet,m'_\bullet) \,\delta(m) \; + & &\\
\nonumber
\Big(\sum_{i\in \ST-\ST'} m_i \big(r \,h^0(E_i(m))-r_i
P(m)-sr_i\delta(m) \big) \Big) & \leq& \\
\nonumber
\Big(\sum_{i\in \ST'} m_i \big(r \,P_{E_i}(m)-r_i P(m) \big) \Big) +
\mu(\varphi,E'_\bullet,m'_\bullet) \,\delta(m) & (\leq)& 0 
\end{eqnarray}
The condition that $E_i$ is saturated can be dropped, since
$h^0(E_i(m))\leq h^0(\overline{E}_i(m))$ and
$\mu(\varphi,E_\bullet,m_\bullet)= 
\mu(\varphi,\overline{E}_\bullet,m_\bullet)$, where $\overline{E}_i$
is the saturated subsheaf generated by $\overline{E}_i$ in $E$.

\medskip
\noindent\textbf{3. $\Rightarrow$ 1.}
Let $(E,\varphi,u)\in \cS''_m$ and consider a 
saturated weighted filtration 
$(E_\bullet,m_\bullet)$. Since $E$ is $N_0$-regular, $P(m)=h^0(E(m))$.
If $E_i\in \cS_0$, then $P_{E_i}(m)=h^0(E_i(m))$. Hence hypothesis
3 applied to the subfiltration $(E'_\bullet,m'_\bullet)$ obtained
by those terms such that $E_i\in \cS_0$ implies
$$
\Big( \sum_{E_i\in \cS_0} m_i\big( r^i P(m) -r P_{E^i}(m) \big) \Big) +
\mu(\varphi,E'_\bullet,m'_\bullet) \,\delta(m) \;(\leq)\; 0\,.
$$
This is equivalent to
$$
\Big(\sum_{E_i\in \cS_0} m_i\big( r P_{E_i}(m) -r_i P(m) \big) \Big)+
\mu(\varphi,E'_\bullet,m'_\bullet) \,\delta(m) \;(\leq)\; 0\, ,
$$
and by lemma \ref{lem2.8}, this is in turn equivalent to
\begin{equation}
\label{2.1.6}
\Big(\sum_{E_i\in \cS_0} m_i\big( r P_{E_i} -r_i P \big) \Big)+
\mu(\varphi,E'_\bullet,m'_\bullet)\, \delta \;(\preceq)\; 0\, .
\end{equation}

If $E_i\notin \cS_0$, then the second alternative of lemma
\ref{lem2.4} holds, and then 
\begin{equation}
\label{2.1.7}
r P_{E_i} -r_i P+s \,r_i\, \delta  \;\prec\; 0\, .
\end{equation}
Using lemma \ref{estimate},  (\ref{2.1.6}) and (\ref{2.1.7})
$$
\Big( \sum_{i=1}^t m_i\big( r P_{E_i} -r_i P \big) \Big)+
\mu(\varphi,E_\bullet,m_\bullet) \,\delta \;(\preceq)\; 0\, .
$$
Again, we can drop the condition that the filtration is saturated,
and this finishes the proof of theorem \ref{thm2.1}

\end{proof}

\begin{corollary}
\label{col2.9}
Let $(E,\varphi,u)$ be $\delta$-semistable, $m\geq N_0$, and assume that
there is a weighted filtration $(E_\bullet,m_\bullet)$ with
$$
\Big(\sum_{i=1}^t m_i\big( r h^0(E_i(m)) -r_i P(m) \big)\Big) +
\mu(\varphi,E_\bullet,m_\bullet) \,\delta(m) \;=\; 0\, .
$$
Then $E_i\in \cS_0$ and $h^0(E_i(m))=P_{E_i}(m)$ for all $i$.
\end{corollary}

\begin{proof}
By the proof of the part $(1\Rightarrow 2)$ of theorem \ref{thm2.1},
if we have this equality then all inequalities in (\ref{chain})
are equalities, hence $\ST'=\ST$, $E_i\in \cS_0$ for all $i$, and the result
follows.
\end{proof}

Note that in theorem \ref{thm2.1} we are assuming that $E$ is torsion
free. To handle the general case, we will use the following lemma

\begin{lemma}
\label{lem1.11}
Fix $u\in \param$.
Let $(E,\varphi,u)$ be a \tensor. 
Assume that there is a family $(E_t,\varphi_t,u)_{t\in C}$
parametrized by a smooth curve $C$ such
that $(E^{}_0,\varphi^{}_0,u)=(E,\varphi,u)$ and $E_t$ is torsion free for
$t\neq 0$. Then there exists a \tensor\ $(F,\psi,u)$, a homomorphism 
$$
(E,\varphi,u)\too (F,\psi,u)
$$
such that $F$ is
torsion free with $P_E=P_F$, and an exact sequence
$$
0 \too T(E) \too E \stackrel{\beta}{\too} F,
$$
where $T(E)$ is the torsion subsheaf of $E$.
\end{lemma}

\begin{proof}
The family is given by a tuple $(E_C,\varphi_C,u_C,N)$
as in (\ref{family}), where $u_C$ is the constant map from
$C$ to $\param$ with constant value $u$. 
Shrinking $C$, we can assume that $N$ is trivial.
Let $U=(X\times C)-\operatorname{Supp}(T(E_0))$. 
Let $F_C=j_*(E_C|_U)$. Since it has no $C$-torsion, 
$F_C$ is flat
over $C$. The natural map $\wt{\beta}:E_C \to F_C$ 
is an isomorphism on $U$,
hence we have a homomorphism $\psi_U:=\varphi_C|_U$ on $U$, 
and this extends to a
homomorphism $\psi_C$ on $X\times C$ 
because $\overline{u_C}^* \sfixed$
is locally free. Finally define
$(F,\psi)=(F_0,\psi_0)$, and let $\beta$ be the homomorphism
induced by $\wt{\beta}$.
\end{proof}

\section{GIT construction}
\label{secgitconstruction}

Let $N\geq N_0$ be large enough so that for all $i>0$, 
all line bundles $L$ of degree $d$, 
all locally free sheaves $\fixed_u$ in the family parametrized by
$\param$, 
and all $m>N$, 
we have $h^i(L^{\otimes b}\otimes
\fixed_u (sm))=0$ and $L^{\otimes b}\otimes
\fixed_u (sm)$ is generated by global sections.

Fix $m\geq N$ and let $V$ be a vector 
space of dimension $p=P(m)$. 
The choice of $m$ implies that if $(E,\varphi,u)$ is 
$\delta$-semistable, then $E(m)$ is generated by global sections
and $h^i(E(m))=0$ for $i>0$.
Let $(g,E,\varphi,u)$ be a tuple where $(E,\varphi,u)$ is a 
$\delta$-semistable \tensor\ and $g$ is an isomorphism 
$g:V\to H^0(E(m))$. This induces a quotient
\begin{equation}
\label{cons.1}
q: V\otimes \SO_X(-m)\surj E \, .
\end{equation}
Let $\SH$ be the Hilbert scheme of quotients of $V\otimes \SO_X(-m)$
with Hilbert polynomial $P$. Let $l>m$ be an integer, 
and $W=H^0(\SO_X(l-m))$. The  quotient $q$ 
induces homomorphisms
$$
\begin{array}{rccl}
q\,\,: & V\otimes \SO_X(l-m)&\surj& E(l) \\
q'\,:& V\otimes W&\to& H^0(E(l)) \\
q'':& \bigwedge{}^{P(l)}(V\otimes W) &\to& 
\bigwedge{}^{P(l)} H^0(E(l))\;\cong \; \CC
\end{array}
$$
If $l$ is large enough, these homomorphisms are surjective, and
give Grothendieck's embedding
$$
\SH \;\too\; \PP(\bigwedge{}^{P(l)}(V^\vee\otimes W^\vee)),
$$
and hence a very ample line bundle $\SO_\SH(1)$ on $\SH$ (depending on
$m$ and $l$).

The tuple $(g,E,\varphi,u)$ induces a linear map
\begin{equation}
\label{cons.2}
\Phi: (V^{\otimes s})^{\oplus c} \;\too\; 
H^0((E(m)^{\otimes s})^{\oplus c}) \;\too\;
H^0((\det E)^{\otimes b} \otimes \fixed_u(sm)).
\end{equation}
Fix a Poincare bundle $\SP$ on $J\times X$,
where $J=\Pic^d(X)$. Fix an isomorphism 
$$
\beta: \det(E) \too \SP|_{\{\det(E)\}\times X}.
$$
Then $\Phi$ induces a quotient
$$
(V^{\otimes s})^{\oplus c} \otimes H^0(\SP|_{\{\det(E)\}\times X}
^{\otimes b} \otimes D_u(sm))^\vee \too \CC.
$$
Choosing a different isomorphism $\beta$ will only change 
this quotient by a scalar, so we get a well defined point
$[\Phi]$ in $P$, where $P$ is the projective
bundle over $J\times \param$ defined as
$$
P=\PP\Big( \big( (V^{\otimes s})^{\oplus c}\big)^\vee \otimes
\pi_{J\times \param *}^{}\big(
\pi_{X \times J}^*\SP^{\otimes b} \otimes 
\pi_{X \times \param}^* \sfixed(sm) \big) \Big) 
\;{\too}\; J\times \param \; .
$$
where $\pi^{}_{X\times J}$ (respectively $\pi^{}_{J\times \param}$,...) 
denotes the natural projection
from $X\times J\times \param$ to $X\times J$ 
(respectively $J\times \param$,...).
Recall that $\sfixed(m):= \sfixed \otimes \pi^*_X\SO_X(m)$.
Note that 
$\pi^{}_{J\times \param *}\big(
\pi_{X\times J}^*\SP^{\otimes b} \otimes 
\pi_{X\times \param}^* \sfixed (sm) \big)$
is locally free because of the choice of $m$.
Replacing $\SP$ with
another Poincare bundle defined by tensoring with the pullback of a
sufficiently positive line bundle on $J$, we can assume that
$\SO_P(1)$ is very ample (this line bundle depends on $m$).

A point $(q,[\Phi])\in \SH\times P$ associated to a tuple 
$(g,E,\varphi,u)$
has the property that the homomorphism $\Phi$ in (\ref{cons.2}) 
composed with evaluation factors as
\begin{equation}
\label{cons.3}
\xymatrix{
{(V^{\otimes s})^{\oplus c}\otimes \SO_X} \ar[d]^{\Phi\qquad} 
\ar@{->>}[r]^{(q^{\otimes s})^{\oplus c}} &
{(E(m)^{\otimes s})^{\oplus c} \ar@(d,r)[ddl]^{\varphi}} \\ 
{H^0((\det E)^{\otimes b} \otimes \fixed_u(sm))\otimes \SO_X} 
\ar[d]^{\ev} & \\
{(\det E)^{\otimes b} \otimes \fixed_u(sm)}}
\end{equation}

Consider the relative version of the homomorphisms in
(\ref{cons.3}), i.e. the commutative diagram on $X\times\SH\times P$
\begin{equation}
\label{cons.3rel}
\xymatrix{
{0} \ar[r] &
{\SK} \ar[r] \ar[rd]_{f}&
{(V^{\otimes s})^{\oplus c}\otimes \SO^{}_{X\times\SH\times P}} 
\ar[r] \ar[d]^{\Phi_{\SH\times P}} &
{(p_{X\times \SH}^*E_{\SH}^{}(m)^{\otimes s})^{\oplus c}} \ar[r] &
{0}\\
 & & {\text{\makebox[0.5cm][l]{$\SA \;:=\; p^*_{X\times J} \SP^{\otimes
b}\otimes p_{X\times \param}^* \sfixed  \otimes
p_{X}^*\SO^{}_X(sm)$}}}
}
\end{equation}
where $p^{}_{X\times J}$ (respectively $p^{}_X$, ...) 
denotes the natural projection
from $X\times\SH\times P$ to $X\times J$ (respectively $X$, ...),
$E_\SH$ is the tautological sheaf on $X\times \SH$, and
$\Phi_{\SH\times P}$ is the relative version of the 
composition $\ev\circ\Phi$ in diagram (\ref{cons.3}).

The points $(q,[\Phi])$ where the restriction
$\Phi_{\SH\times P}|_{X\times (q,[\Phi])}$ factors through
$(E(m)^{\otimes s})^{\oplus c}$ (as in \ref{cons.3})
are the points where $f|_{X\times (q,[\Phi])}$
is identically zero. 
We will need the following

\begin{lemma}
\label{zeroes}
Let $Y$ be a scheme, and let $f:\SG\to \SF$ be a 
homomorphism of coherent sheaves on $X\times Y$. Assume that
$\SF$ is flat over $Y$. 
Then there is a unique closed subscheme $Z$
satisfying the following universal property:
given a Cartesian diagram
$$
\xymatrix{
{X\times S} \ar[r]^-{\overline{h}} \ar[d]^{{p_S}}& 
{X\times Y} \ar[d]^{p}
\\
{S} \ar[r]^-{h} & {Y} 
}
$$
$\overline{h}^* f=0$ if and only if $h$ factors through $Z$.

\end{lemma}

\begin{proof}
Uniqueness is clear.
To show existence, assume that $\SO_X(1)$ is very ample 
(taking a multiple if necessary) and
let $p:X\times Y\to Y$ be the projection to the second factor.
Since $\SF$ is $Y$-flat, taking $m'$ 
large enough, $p_*^{} \SF(m')$ is locally 
free (recall $\SF(m')= \SF\otimes p_X^*\SO_X^{}(m')$). 
The question is local on $Y$,
so we can
assume, shrinking $Y$ if necessary, 
that $Y=\Spec A$ and $p_*^{} \SF(m')$ is given by 
a free $A$-module.
Now, since $Y$ is affine, the homomorphism
$$
p^{}_*f(m'): p_*^{} \SG(m') \too p_*^{} \SF(m')
$$
of sheaves on $Y$ is equivalent to a homomorphism of 
$A$-modules
$$
M \stackrel{(f_1,\ldots,f_n)}{\too} A\oplus \cdots\oplus A
$$
The zero locus of $f_i$ is defined by the ideal $I_i\subset A$
image of $f_i$, thus the zero scheme of $(f_1,\ldots,f_n)$ 
is given by the ideal $I=\sum I_i$,
hence $Z'_{m'}$ is a closed subscheme.

Since $\SO_X(1)$ is very ample, if $m''>m'$
we have an injection $p_*^{} \SF(m')\inj p_*^{} \SF(m'')$
(and analogously for $\SG$), hence $Z_{m''}\subset Z_{m'}$,
and since $Y$ is noetherian, there exists $N'$ such
that, if $m'>N'$,
we get a scheme $Z$ independent of $m'$.

To check the universal property first we will show that if
$\overline{h}^*f=0$ then $h$ factors through $Z$.  Since the question
is local on $S$, we can take $S=\Spec(B)$, $Y=\Spec(A)$, 
and the morphism $h$ is locally
given by a ring homomorphism $A\to B$.  Since $\SF$ is flat over $Y$,
for $m'$ large enough the natural homomorphism $\alpha:h^*p_*\SF(m')\to
{p_S}_*\overline{h}^*\SF(m')$ (defined as in \cite[Th. III 9.3.1]{Ha})
is an isomorphism.  Indeed, for $m'$ sufficiently large,
$H^i(X,\SF_y(m'))=0$ and $H^i(X,\overline{h}^* (\SF(m'))_s)=0$ for all
points $y\in Y$, $s\in S$ and $i>0$, and since $\SF$ is flat, this
implies that $h^*p_*\SF(m')$ and ${p_S}_*\overline{h}^*\SF(m')$ are
locally free. Then, to prove that the homomorphism $\alpha$ 
is an isomorphism,
it is enough to check it at the fiber of every $s\in S$, but this
follows from \cite[Th. III 12.11]{Ha} or \cite[II \S 5 Cor. 3]{Mu}.

Hence the commutativity of
the diagram
$$
\xymatrix{
{{p_S}_*\overline{h}^*\SG(m')} 
\ar[rrr]^{{p_S}_*\overline{h}^*f(m')=0}&&&
{{p_S}_*\overline{h}^*\SF(m')}  \\
{h^*p_*\SG(m')} \ar[rrr]^{h^*p_*f(m')} \ar[u] &&&
{h^*p_*\SF(m')} \ar[u]_{\isom}
}
$$
implies that $h^*p_*f(m')=0$. This means that 
for all $i$,  in the
diagram
$$
\xymatrix{
{M} \ar[r]^{f_i}\ar[d] & {A} \ar[r]\ar[d] & {A/I_i}\ar[d] \ar[r] & {0}\\
{M\otimes_A B} \ar[r]^-{f_i\otimes B} 
& {B} \ar[r] & {A/I_i\otimes_A B}\ar[r] &{0}
}
$$
it is $f_i\otimes B=0$.
Hence the image $I_i$ of
$f_i$ is in the kernel $J$ of $A\to B$. Therefore $I\subset J$,
hence $A\to B$ factors through $A\to A/I$, which means that
$h:S\to Y$ factors through $Z$.

Now we show that if we take $S=Z$ and $h:Z\inj Y$ the inclusion,
then $\overline{h}^* f=0$.
By definition of $Z$ we have $h^*p_* f(m')=0$ for any $m'$ with 
$m'>N'$.
Showing that $\overline{h}^* f=0$
is equivalent to showing that
$$
\overline{h}^* f(m')  : \overline{h}^* \SG(m') \too \overline{h}^* \SF(m')
$$ 
is zero for some $m'$.
Take $m'$ large enough so that 
$\ev:p^*p_*\SG(m')\to \SG(m')$ is surjective.
By the right exactness of $\overline{h}^*$ the homomorphism 
$\overline{h}^*\ev$ is still surjective. The commutative
diagram
$$
\xymatrix{
{\overline{h}^* \SG(m')} \ar[rrr]^{\overline{h}^* f(m')} &&& 
{\overline{h}^* \SF(m')}\\ 
{\overline{h}^*p^*_{}p_*  \SG(m')} 
\ar[rrr]^{\overline{h}^*p^*_{}p_*  f(m')}
\ar@{>>}[u]^{\overline{h}^*\ev} &&& 
{\overline{h}^*p^*_{}p_* \SF(m')} \ar[u]\\
{p_S^{*}h^*_{}p_*  \SG(m')} 
\ar[rrr]^{p_S^{*}h^*_{}p_*  f(m')=0} \ar@{=}[u]&&& 
{p_S^{*}h^*_{}p_* \SF(m')} \ar@{=}[u]
}
$$
implies $\overline{h}^*f(m')=0$, hence $\overline{h}^*f=0$.

\end{proof}

Let $Z'$ be the scheme given by this lemma for
$Y=\SH\times P$ and
the homomorphism $f:\SK\to \SA$. Then
$\overline{i}^*f=0$, and there is a 
commutative diagram on $X \times Z'$
\begin{equation}
\label{cons.8}
\xymatrix{
{\overline{i}^* \SK} \ar[r] \ar[rd]_{\overline{i}^* f}&
{(V^{\otimes s})^{\oplus c}\otimes \SO^{}_{X\times Z'}} 
\ar[r] \ar[d]^{\overline{i}^*\widetilde\Phi} &
{(\overline{i}^* p_{X\times \SH}^*E^{}_\SH(m)^{\otimes s})^{\oplus c}} \ar[r] 
\ar[ld]^{\widetilde\varphi} &
{0}\\
  &  {\overline{i}^*\SA}
}
\end{equation}
and hence there is a universal family of based \tensors\ parametrized
by $Z'$
\begin{equation}
\label{3.5bis}
\varphi_{Z'}:E^{}_{Z'}{}^{\otimes s} 
\too (\det E_{Z'})^{\otimes b}\otimes p^*_{Z'}N
\end{equation}

Given a point $(q,[\Phi])$ in $Z'$, 
using the tautological family (\ref{3.5bis})
we can recover the tuple $(q,E,\varphi,u)$ up to isomorphism,
and if $H^0(q(m)):V\to H^0(E(m))$ is an isomorphism,
then we recover the tuple $(g=H^0(q(m)),E,\varphi,u)$ up to isomorphism,
i.e. if $(g',E',\varphi',u')$ is another tuple corresponding to
the same point, then there exists an isomorphism $(f,\alpha)$
between $(E,\varphi,u)$ and $(E',\varphi',u')$ as in
(\ref{isomorphism}), and $H^0(f(m))\circ q=q'$.

Let $Z\subset Z'$ be the
closure of the points associated to $\delta$-semistable \tensors.
Let $p^{}_\SH$ and $p^{}_P$ be the projections of $Z$ to $\SH$ and $P$, and
define a polarization on $Z$ by 
\begin{equation}
\label{son1n2}
\SO_Z(n_1,n_2):=p^*_\SH \SO^{}_\SH(n_1) \otimes 
p^*_P \SO^{}_P(n_2),
\end{equation}
where $n_1$ and $n_2$ are integers with
\begin{equation}
\label{n2n1}
\frac{n_2}{n_1}=\frac{P(l)\delta(m)-\delta(l)P(m)}{P(m)-s\delta(m)}
\end{equation}
The projective scheme $Z$ is preserved by the natural $\slv$ action,
and this action has a natural linearization on $\SO_Z(n_1,n_2)$, using
the natural linearizations on $\SO_\SH(1)$ and $\SO_P(1)$. 

We have seen that the points of $Z$ for which $H^0(q(m))$ is an
isomorphism correspond
(up to isomorphism) to tuples $(g,E,\varphi,u)$, where $g$ is an
isomorphism between $V$ and $H^0(E(m))$.
To get rid of the choice of $g$, we have to take the quotient 
by $\glv$, but if $\lambda\in \CC^*$, $(g,E,\varphi,u)$ and 
$(\lambda g,E,\varphi,u)$ correspond to the same point, 
and hence it is enough to divide by the action of $\slv$.
In fact, the moduli
space will be the GIT quotient of $Z$ by $\slv$.

In proposition \ref{prop3.1}, we will identify the 
GIT-(semi)stable points  in $Z$ using
the Hilbert-Mumford criterion. In theorem \ref{thm3.2} we relate
filtrations of sheaves with filtrations of the vector space $V$
to prove that GIT-(semi)stable points
of $Z$ coincide with the points associated to $\delta$-(semi)stable
points.

A nonconstant group homomorphism $\lambda:\CC^*\too \slv$ is called a 
one-parameter subgroup
of $\slv$. If $\slv$ acts on a projective scheme $Y$ with a given
linearization, we denote by $\mu(y,\lambda)$ the minimum weight of the
action of $\lambda$ on $y\in Y$.

A weighted filtration $(V_\bullet,m_\bullet)$ 
of the vector space $V$ is a filtration of
vector spaces
\begin{equation}
\label{filtV}
0 \subsetneq V_1 \subset V_2 \subset \;\cdots\; \subset V_t \subsetneq
V_{t+1}=V,
\end{equation}
and positive numbers $m_1,\, m_2,\ldots , \,m_t>0$. If $t=1$ (one step
filtration), then we will take $m_1=1$.
Consider the vector of $\CC^p$ defined as $\Gamma=\sum_{i=1}^{t} 
m_i \Gamma^{(\dim(V_i))}$, where
$$
\Gamma^{(k)}=\big( \overbrace{k-p,\dots,k-p}^k,
\overbrace{k,\dots,k}^{p-k} \big)
\qquad
(1\leq k < p)\,.
$$
Now let $\SI=\{1,\ldots,t+1\}^{\times s}$ be the set of all
multi-indexes $I=(i_1,\ldots,i_s)$, and define
\begin{equation}
\label{V1}
\mu(\Phi,V_\bullet,m_\bullet) \;=\; \min_{I\in \SI} \big\{
\Gamma_{\dim V_{i_1}}+\dots+\Gamma_{\dim V_{i_s}}: \,
\Phi|_{(V_{i_1}\otimes\cdots \otimes V_{i_s})^{\oplus c}}\neq 0 \;
  \big\} 
\end{equation}
As we did in the proof of lemma \ref{estimate},
if $I$ is
the multi-index giving minimum in (\ref{V1}), we will denote
by $\epsilon_i(\Phi,V_\bullet,m_\bullet)$ 
(or just $\epsilon_i(V_\bullet)$ if the rest of the data is
clear from the context)
the number of 
elements $k$ of the multi-index $I$ such
that $\dim V_k\leq \dim V_i$. Then we have, as in (\ref{epsilonE})
\begin{equation}
\label{epsilonV}
\mu(\Phi,V_\bullet,m_\bullet) \;=\;
\sum_{i=1}^t m_i \big( s\dim V_i-\epsilon_i(V_\bullet)\dim V \big).
\end{equation}

Given a subspace $V'\subset V$ and a quotient
$q:V\otimes\SO_X(-m)\surj E$, we define the subsheaf $E_{V'}$ of $E$ as
the image of the restriction of $q$ to $V'$
$$
\xymatrix{
{V\otimes \SO_X(-m)} \ar@{>>}[r] & {E} \\
{V'\otimes \SO_X(-m)} \ar@{>>}[r] \ar@{^{(}->}[u] & 
{E_{V'}} \ar@{^{(}->}[u]\\
}
$$
In particular, $E_{V'}(m)$ is generated by global sections.

On the other hand, if the quotient $q:V\otimes\SO_X(-m)\surj E$ 
induces an injection $V\inj H^0(E(m))$ (we will later show that all 
quotients
coming from GIT-semistable points of $Z$ satisfy this property),
and if $E'\subset E$ is a subsheaf, we define 
$$
V_{E'}=V \cap H^0(E'(m)).
$$
The following two lemmas are easy to check

\begin{lemma}
\label{A1}
Given a point $(q,[\Phi])\in Z$ such that $q$ induces an injection
$V\inj H^0(E(m))$, and a weighted filtration $(E_\bullet,m_\bullet)$ 
of $E$, we have:
\begin{enumerate}

\item $E_{V_{E_i}} \subset E_i$

\item If 
$\varphi|_{(E_{i_1}\otimes \dots \otimes E_{i_s})^{\oplus c}}=0$,
then 
$\Phi|_{(V_{E_{i_1}}\otimes \dots \otimes V_{E_{i_s}})^{\oplus c}}=0$

\item $\sum_{i=1}^{t}-m_i \, \epsilon_i(\varphi,E_\bullet,m_\bullet) 
\leq \sum_{i=1}^{t}-m_i \, \epsilon_i(\Phi,V_{E_\bullet},m_\bullet)$
\end{enumerate}
Furthermore, if $q$ induces an isomorphism $V\cong H^0(E(m))$,
all $E_i$ are $m$-regular and all 
$E_{i_1}\otimes \cdots \otimes E_{i_s}$
are $sm$-regular, then (1) becomes an equality,
(2) becomes ``if and only if'' and (3) 
an equality.
\end{lemma}

\begin{lemma}
\label{A2}
Given a point $(q,[\Phi])\in Z$ such that $q$ induces an injection
$V\inj H^0(E(m))$, and a weighted filtration $(V_\bullet,m_\bullet)$ 
of $V$, we have:
\begin{enumerate}

\item $V_i\subset V_{E_{V_i}}$

\item $\varphi|_{(E_{V_{i_1}}\otimes \dots \otimes E_{V_{i_s}})^{\oplus
c}}=0$ if and only if
$\Phi|_{(V_{i_1}\otimes \dots \otimes V_{i_s})^{\oplus c}}=0$

\item $\sum_{i=1}^{t}-m_i \,\epsilon_i(\varphi,E_{V_\bullet},m_\bullet) 
= \sum_{i=1}^{t}-m_i \,\epsilon_i(\Phi,V_\bullet,m_\bullet)$
\end{enumerate}
\end{lemma}

\begin{proposition}
\label{prop3.1}
For sufficiently large $l$, the point $(q,[\Phi])$ in $Z$ is
GIT-(semi)stable with respect to $\SO_Z(n_1,n_2)$ if and only if for
every weighted filtration $(V_\bullet,m_\bullet)$ of $V$
\begin{equation}
\label{eqprop3.1}
n_1\Big(\sum_{i=1}^{t}m_i \big( \dim V_i P(l)- \dim V \,P_{E_{V_i}}(l)
\big)\Big) \,+\, 
n_2\, \mu(\Phi,V_\bullet,m_\bullet) \;(\leq)\; 0.
\end{equation}
Furthermore, there is an integer $A_2$ (depending only on $m$, $P$,
$s$, $b$, $c$ and $\sfixed$) such that it is enough to consider weighted
filtrations with $m_i\leq A_2$.
\end{proposition}

\begin{proof}
Given $m$, the sheaves $E_{V'}$ for $V'\subset V$ form a bounded
family, so if $l$ is large enough, we will have
$$
\dim q'(V\otimes W)\,=\,h^0(E_{V'}(l))\,=\,P_{E_{V'}}(l)
$$
for all subspaces $V'\subset V$.
By the Hilbert-Mumford criterion, a point is GIT-(semi)stable if and
only if for all one-parameter subgroups $\lambda$ of $\slv$,
$$
\mu((q,[\Phi]),\lambda)=\,n_1\mu(q,\lambda) +n_2\,\mu([\Phi],\lambda)
\;(\leq)\;0.
$$
A one-parameter subgroup of $\slv$ is equivalent to a basis 
$\{e_1,\ldots,e_p\}$ of $V$
and a vector $\Gamma=(\Gamma_1,\ldots,\Gamma_p)\in \CC^p$ with
$\Gamma_1\leq \ldots \leq \Gamma_p$. This defines a 
weighted filtration $(V_\bullet,m_\bullet)$ of $V$ as follows:
let $\lambda_1<\ldots<\lambda_{t+1}$ be the different values of
$\Gamma_k$, let $V_i$ be the vector space generated by all $e_k$ such
that $\Gamma_k\leq \lambda_i$, and let
$m_i=(\lambda_{i+1}-\lambda_i)/p$.
Denote $\SI'=\{1,\ldots,t+1\}^{\times P(l)}$. 
We have (\cite{Si} or \cite{H-L2})
\begin{eqnarray*}
\mu(q,\lambda)&=&
\min_{I\in\SI'} 
\big\{ \Gamma_{\dim V_{i_1}}+
\dots +\Gamma_{\dim V_{i_{P(l)}}} : q''|_{(V_{i_1}\otimes W)
\bigwedge \dots \bigwedge (V_{i_{P(l)}}\otimes W)} \neq 0 \big\} \\
& = &
 \sum_{i=1}^{t}m_i \big( \dim V_i\, P(l) - 
\dim V \,P_{E_{V_i}}(l) \big)
\end{eqnarray*}
and
\begin{eqnarray*}
\mu([\Phi],\lambda) &=& \min_{I\in \SI} \big\{
\Gamma_{\dim V_{i_1}}+\dots+\Gamma_{\dim V_{i_s}}: \,
\Phi|_{(V_{i_1}\otimes\cdots \otimes V_{i_s})^{\oplus c}}\neq 0 \big\}\\
 &=& \mu(\Phi,V_\bullet,m_\bullet)
\end{eqnarray*}

The last statement follows from an argument similar to the proof of
lemma \ref{finiteE}, with $\ZZ^r$ replaced by $\ZZ^p$.

\end{proof}

\begin{proposition}
\label{thm3.2a}
The point
$(q,[\Phi])$ is GIT-(semi)stable if
and only if for all weighted filtrations $(E_\bullet,m_\bullet)$ of $E$
\begin{equation}
\label{critform}
\begin{array}{r}
\sum_{i=1}^{t} m_i\Big( \big(\dim V_{E_i}
-\epsilon_i(E_\bullet)\delta(m)\big)
\big(P-s\delta\big) \hspace{2cm}\\
\mbox{}-\big(P_{E_i}-\epsilon_i(E_\bullet)\delta\big)\big(\dim V
-s\delta(m)\big) 
\Big) 
\;(\preceq)\; 0
\end{array}
\end{equation}
Furthermore, if $(q,[\Phi])$ is GIT-semistable, then the induced map
$f_q:V\to H^0(E(m))$ is injective.
\end{proposition}

\begin{proof}
First we prove that if $(q,[\Phi])$ is GIT-semistable, then the
induced linear map $f_q$ is injective. Let $V'$ be its
kernel and consider the filtration $V'\subset V$. We have $E_{V'}=0$
and $\mu(\Phi,V'\subset V)=s\dim V'$. Applying proposition
\ref{prop3.1} we have
$$
n_1 \dim V' P(l) + n_2 s \dim V' \leq 0,
$$
and hence $V'=0$.

Using (\ref{n2n1}) and (\ref{epsilonV}), the inequality of proposition
\ref{prop3.1} becomes
\begin{equation}
\label{3.2.1}
\begin{array}{r}
\sum_{i=1}^{t} m_i\Big( \big(\dim V_i-\epsilon_i(V_\bullet)\,\delta(m)\big)
\big(P(l)-s\delta(l)\big) \hspace{2cm}\\
\mbox{}-\big(P_{E_{V_i}}(l)-\epsilon_i(V_\bullet)\delta(l)\big)
\big(\dim V -s\delta(m)\big) \Big)
\;(\leq)\; 0
\end{array}
\end{equation}

An argument similar to lemma \ref{lem2.8} (using $A_2$ instead of
$A_1$) shows that we can take $l$
large enough 
(depending only on $m,\,s,\,b,\,c,\,P,\,\sfixed$ and $\delta$), 
so that this inequality 
holds for $l$ if and only if it holds as an inequality of 
polynomials.

Now assume that $(q,[\Phi])$ is GIT-(semi)stable. 
Take a weighted filtration $(E_\bullet,m_\bullet)$ of $E$. 
Then lemma \ref{A1} and (\ref{3.2.1}) applied to 
the associated 
weighted filtration $(V_{E_\bullet},m_\bullet)$ of $V$ give 
(\ref{critform}).

On the other hand, assume that (\ref{critform}) holds. 
Take a weighted filtration $(V_\bullet,m_\bullet)$ of $V$. 
Then lemma \ref{A2} and (\ref{critform}) applied to 
the associated 
weighted filtration $(E_{V_\bullet},m_\bullet)$ of $E$ give
(\ref{3.2.1}), and it follows that $(q,[\Phi])$ is GIT-(semi)stable.

\end{proof}

\begin{theorem}
\label{thm3.2}
Assume $m> N$. For $l$ sufficiently large, a point $(q,[\Phi])$ in $Z$ is 
GIT-(semi)stable if and only
if the corresponding \tensor\ $(E,\varphi,u)$ is $\delta$-(semi)stable
and the linear map $f_q:V\too H^0(E(m))$ induced by $q$ is an isomorphism.
\end{theorem}

\begin{proof}
We prove this in two steps:

\medskip
\noindent
\textbf{Step 1.}\textit{
$(q,[\Phi])$  GIT-semistable
$\Longrightarrow$ $(E,\varphi,u)$ $\delta$-semistable and 
$q$ induces an isomorphism $V\isom H^0(E(m))$.
}
\medskip

The leading coefficient of (\ref{critform}) gives
$$
\sum_{i=1}^{t} m_i\Big( \big(\dim V_{E_i} 
-\epsilon_i(E_\bullet)\delta(m)\big)r
-r_i\big(\dim V -s\delta(m)\big) \Big) 
\;\leq\; 0\, .
$$
Note that even if $(q,[\Phi])$ is GIT-stable, here we only get weak
inequality. This implies
\begin{equation}
\label{3.2.4}
\Big(\sum_{i=1}^{t} m_i\big(r^i P(m) -r h^0(E^i(m))\big)\Big)
+\mu(\varphi,E_\bullet,m_\bullet)\delta(m)
\;\leq\; 0\, .
\end{equation}
To be able to apply theorem \ref{thm2.1}, we still need to show that $E$
is torsion free. By lemma \ref{lem1.11}, there exists a \tensor\ 
$(F,\psi,u)$
with $F$ torsion free such that $P_E=P_F$ and an exact sequence
$$
0 \too T(E) \too E \stackrel{\beta}{\too} F\,.
$$
Consider a weighted filtration $(F_\bullet,m_\bullet)$ of $F$. 
Let $F^i=F/F_i$, and
let $E^i$ be the image of $E$ in $F^i$. Let $E_i$ be the kernel
of $E \to E^i$. Then $\rk(F_i)=rk(E_i)=r_i$, 
$h^0(F^i(m))\geq h^0(E^i(m))$, and
$\mu(\psi,F_\bullet,m_\bullet)= \mu(\varphi,E_\bullet,m_\bullet)$. 
Using this and applying
(\ref{3.2.4}) to $E_i$ we get
$$
\Big(\sum_{i=1}^{t} m_i\big( r^i P(m) -r h^0(F^i(m))\big)  \Big)
+\mu(\varphi,F_\bullet,m_\bullet)\,\delta(m)
\;\leq\; 0,
$$
and hence theorem \ref{thm2.1} implies that $(F,\psi,u)$ is 
$\delta$-semistable.

Next we will show that $T(E)=0$, and hence, since $P_E=P_F$,
we will conclude that $(E,\varphi,u)$ is isomorphic to $(F,\psi,u)$.
Define $E''$ to be the image of $E$ in $F$. Then
$$
P(m)-s\delta(m) \,=\,h^0(F(m))-s\delta(m)\,\geq\,
h^0(E''(m))-s\delta(m) \,\geq\,
P(m)-s\delta(m),
$$
where the last inequality follows from (\ref{3.2.4}) applied to
the one step filtration $T(E)\subset E$. Hence equality holds at all places and
$h^0(F(m))=h^0(E''(m))$. Since $F$ is globally generated, $F=E''$, and
hence $T(E)=0$.

Finally, we have seen that $f_q$ is injective, and since $(E,\varphi)$
is $\delta$-semistable, $\dim V = h^0(E(m))$, hence $f_q$ is an
isomorphism.

\medskip
\noindent
\textbf{Step 2.}\textit{
$(E,\varphi,u)$ $\delta$-stable
(respectively strictly $\delta$-semistable) and $q$ induces an isomorphism 
$f_q:V\isom H^0(E(m))$
$\Longrightarrow$ $(q,[\Phi])$  GIT-stable (respectively strictly semistable).
}
\medskip

Since $f_q$ is an isomorphism, we have $V_{E'}=H^0(E'(m))$ for any subsheaf
$E'\subset E$. Then theorem \ref{thm2.1} implies that for all weighted
filtrations
\begin{equation}
\label{3.2.5}
\Big(\sum_{i=1}^t m_i \big(r \dim V_{E_i}-r_i P(m)\big) \Big) +
\mu(\varphi,E_\bullet,m_\bullet) \,\delta(m) \;(\leq)\; 0
\end{equation}
If the inequality is strict, then 
$$
\begin{array}{r}
\sum_{i=1}^{t} m_i \Big( \big(\dim V_{E_i} -\epsilon_i(E_\bullet)\delta(m)
  \big) \big(P-s\delta\big) \hspace{2cm}\\
\mbox{}-\big( P_{E_i}-\epsilon_i(E_\bullet)\delta\big)
  \big(\dim V -s\delta(m)\big)
  \Big) \;\prec\; 0 \,.
\end{array}
$$
If $(E,\varphi,u)$ is strictly $\delta$-semistable, by theorem
\ref{thm2.1} there is a filtration giving equality in (\ref{3.2.5}),
then corollary \ref{col2.9} implies that $h^0(E_i(m))=P_{E_i}(m)$, 
and by lemma \ref{lem2.8}
$$
\Big(\sum_{i=1}^t m_i(r P_{E_i}-r_i P) \Big)+
\mu(\varphi,E_\bullet,m_\bullet)\, \delta \;=\; 0\, ,
$$
and a short calculation using this and (\ref{3.2.5}) 
gives
$$
\begin{array}{r}
\sum_{i=1}^{t} m_i \Big( \big(\dim V_{E_i} -\epsilon_i(E_\bullet)\delta(m)
  \big) \big(P-s\delta\big) \hspace{2cm}\\
\mbox{}-\big( P_{E_i}-\epsilon_i(E_\bullet)\delta\big)
  \big(\dim V -s\delta(m)\big)
  \Big) \;=\; 0 \,.
\end{array}
$$
So we finish by using proposition \ref{thm3.2a}.

\end{proof}

Given a one-parameter subgroup of $\slv$, choose
a basis $\{e_j\}$ of $V$ where it has a diagonal
form 
$$
{\rm diag}(\overbrace{\lambda_1,\ldots,\lambda_1}^{a_1},
\overbrace{\lambda_2,\ldots,\lambda_2}^{a_2-a_1},\ldots,
\overbrace{\lambda_{t+1},\ldots,\lambda_{t+1}}^{a_{t+1}-a_t})
$$
This gives a weighted filtration $(V_\bullet,m_\bullet)$ of $V=H^0(E(m))$
(where $V_i$ is the linear span of $\{e_j\}$ with $j\leq a_i$,
and $m_i=(\lambda_{i+1}-\lambda_i)/p$) and a
splitting $V=\oplus V^i$ of this filtration 
(with $V^i$ the linear span of $\{e_j\}$ with $a_{i-1}<j\leq a_i$).
Defining $E_{V_i}=q(V_i\otimes \SO_X(-m))$ we obtain a 
weigted filtration $(E_\bullet,m_\bullet)$ of $E$.

Now let $(E_\bullet,m_\bullet)$ be a weigted filtration of $E$
and $V=\oplus V^i$ a splitting of the filtration
$V_i=H^0(E(m))$. This gives a one-parameter subgroup
$\lambda$ of $\slv$, defined as $v^i\mapsto t^{\lambda_i}v^i$
for $v^i\in V^i$, with $\lambda_i$ such that 
$m_i=(\lambda_{i+1}-\lambda_i)/p$ and $\sum \lambda_i \dim V^i=0$.

The following proposition will be used to prove the 
criterion for S-equivalence.

\begin{proposition}
\label{bijection}
Assume $m>N$. 
Let $(E,\varphi,u)$ be a $\delta$-semistable \tensor,
$f:V\isom H^0(E(m))$ an isomorphism, and let
Let $(q,[\Phi])\in Z$ be the corresponding 
GIT-semistable point.
The above construction gives a bijection between
one-parameter subgroups of $\slv$ 
with $\mu((q,[\Phi]),\lambda)=0$ on the one hand,
and weighted filtrations $(E_\bullet,m_\bullet)$
of $E$ with 
$$
\Big(\sum_{i=1}^t m_i(r P_{E_{i}}-r_i P) \Big)+
\mu(\varphi,E_\bullet,m_\bullet)\, \delta \;=\; 0
$$
together with a splitting of the filtration $H^0(E_\bullet(m))$
of $V=H^0(E(m))$ on the other hand.
\end{proposition}

\begin{proof}
Let $\lambda$ be a one-parameter subgroup of $\slv$
with
$\mu((q,[\Phi]),\lambda)=0$.
The proof of proposition \ref{prop3.1} then gives 
equality in (\ref{eqprop3.1}).
Using (\ref{n2n1}) (relationship between $n_2/n_1$ and $\delta$), 
(\ref{epsilonV}) (relationship between $\epsilon(V_\bullet)$
and $\mu(V_\bullet$)) and lemma \ref{A2}
(relationship between $\epsilon(V_\bullet)$ and $\epsilon(E_{V_\bullet})$), 
this equality becomes
$$
\begin{array}{r}
\sum_{i=1}^{t} m_i\Big( \big(\dim V_i-\epsilon_i(E_{V_\bullet})\,\delta(m)\big)
\big(P(l)-s\delta(l)\big) \hspace{2cm}\\
\mbox{}-\big(P_{E_{V_i}}(l)-\epsilon_i(E_{V_\bullet})\delta(l)\big)
\big(\dim V -s\delta(m)\big) \Big)
\;=\; 0
\end{array}
$$
We have chosen $l$ so large that this holds if and only if
it holds as a polynomial in $l$, hence taking the leading
coefficient in $l$ we obtain
$$
\sum_{i=1}^{t} m_i\Big( \big(\dim V_{i} 
-\epsilon_i(E_{V_\bullet})\delta(m)\big)r
-r_i\big(\dim V -s\delta(m)\big) \Big) 
\;=\; 0\, .
$$
where $r_i=\rk E_{V_i}$ and $r=\rk E$.
Using (\ref{epsilonE}), this is
$$
\Big(\sum_{i=1}^t m_i \big(r \dim V_{i}-r_i P(m)\big) \Big) +
\mu(\varphi,E_{V_\bullet},m_\bullet) \,\delta(m) \;=\; 0
$$
By lemma \ref{A2}, $V_i\subset V_{E_{V_i}}=H^0(E_{V_i}(m))$, hence
$$
\Big(\sum_{i=1}^t m_i \big(r h^0(E_{V_i}(m))   -r_i P(m)\big) \Big) +
\mu(\varphi,E_{V_\bullet},m_\bullet) \,\delta(m) \;\geq \; 0
$$
but by theorem \ref{thm2.1} this must be nonpositive,
hence $V_i = H^0(E_{V_i}(m))=V_{E_{V_i}}$, and the last inequality
is an equality. By corollary \ref{col2.9},
$E_i\in \cS_0$, and hence 
$h^0(E_{V_i}(m))=P(E_{V_i}(m))$ for
all $i$, and then lemma \ref{lem2.8} gives
$$
\Big(\sum_{i=1}^t m_i \big(r P_{E_{V_i}}   -r_i P\big) \Big) +
\mu(\varphi,E_{V_\bullet},m_\bullet) \,\delta \;= \; 0
$$

Conversely, let $(E_\bullet,m_\bullet)$ be a filtration
with
\begin{equation}
\label{hypo}
\Big(\sum_{i=1}^t m_i \big(r P_{E_{i}}   -r_i P\big) \Big) +
\mu(\varphi,E_{\bullet},m_\bullet) \,\delta \;= \; 0
\end{equation}
together with a splitting of the filtration 
$H^0(E_i(m))$ of $V \cong H^0(E(m))$, and let $\lambda$ be the associated
one-parameter subgroup of $\slv$.
Equation (\ref{hypo}) gives in particular
$$
\Big(\sum_{i=1}^t m_i \big(r P_{E_{i}}(m)   -r_i P(m)\big) \Big) +
\mu(\varphi,E_{\bullet},m_\bullet) \,\delta(m) \;= \; 0
$$
By the proof of implication 3. $\Rightarrow$ 1. in 
theorem \ref{thm2.1}, since we get an equality, it
is $E_i\in \cS_0$ for all $i$, hence $P_{E_i(m)}=h^0(E_i(m))=\dim V_{E_i}$
for all $i$, and the previous equality becomes
\begin{equation}
\label{hypobis}
\Big(\sum_{i=1}^t m_i \big(r \dim V_{E_i}   -r_i \dim V\big) \Big) +
\mu(\varphi,E_{\bullet},m_\bullet) \,\delta(m) \;= \; 0
\end{equation}
Furthermore, the strong version of lemma \ref{A1} gives
$E_i=E_{V_{E_i}}$.
Using (\ref{hypobis}) and (\ref{hypo}), together with
(\ref{epsilonE}) and the strong form of lemma \ref{A1}, we obtain
$$
\sum_{i=1}^{t} m_i\Big( \big(\dim V_{E_i}
-\epsilon_i(V_{E_\bullet})\delta(m)\big)
\big(P-s\delta\big) \hspace{2cm}\\
\mbox{}-\big(P_{E_i}-\epsilon_i(V_{E_\bullet})\delta\big)\big(\dim V
-s\delta(m)\big) 
\Big) =0
$$
Hence, we also get $0$ after evaluating this polynomial
in $l$, but by the proofs of propositions
\ref{prop3.1} and \ref{thm3.2a}, this is
equal to $\mu((q,[\Phi]),\lambda)$.

We have seen that
$V_i=V_{E_{V_i}}$ and $E_i=E_{V_{E_i}}$, and it is
easy to check that this gives a bijection.

\end{proof}

\section{Proof of theorem \ref{mainthm}}
\label{sectheorem}

\begin{proof}[Proof of theorem \ref{mainthm}]
The main ingredient of the proof is theorem
\ref{thm3.2}, showing that GIT-(semi)stable points correspond
to $\delta$-(semi)stable \tensors.
 
Using the notation of section \ref{secgitconstruction},
let $\FM_\delta$ (respectively $\FM^s_\delta$) be the GIT quotient
of $Z$ (respectively $Z^s$) by $\slv$. 
Since $Z$ is projective, $\FM_\delta$ is also projective.
GIT gives that $\FM^s_\delta$ is an open subset of
the projective scheme $\FM_\delta$. The restriction $Z^s\to
\FM^s_\delta$ to the stable part is a geometric quotient, i.e.
the fibers are $\slv$-orbits, and hence the points
of $\FM^s_\delta$ correspond to isomorphism classes of 
$\delta$-stable \tensors.

It only remains to show that $\FM_\delta$ corepresents the functor
${\SM^{}_\delta}$. 
We will follow closely
\cite[Proof of Main Theorem 0.1, p. 315]{H-L2}

Let $(E_T,\varphi_T,u_T,N)$ be a family of
$\delta$-semistable \tensors\ (cf. (\ref{family})) 
parametrized by a scheme $T$.
Then $\SV:=\pi^{}_{T*}(E_T\otimes \pi^*_X \SO^{}_X(m))$ is locally free on
$T$. The family $E_T$ gives a map $\Delta:T\to \Pic^d(X)$, sending
$t\in T$ to $\det E_t$.
Cover $T$ with small open sets $T_i$. For each $i$ we can find an
isomorphism
$$
\beta_{T_i}: \det E_{T_i} \too \overline{\Delta_i}^*\SP
$$
(where $\SP$ is the Poincare bundle in the definition of $P$ at the
beginning of section \ref{secgitconstruction}),
and a trivialization
$$
g_{T_i}: V\otimes \SO_{T_i} \too \SV|_{T_i}.
$$
Using this trivialization we obtain a family of quotients
parametrized by $T_i$
$$
q_{T_i}: V\otimes \pi^*_X\SO^{}_{X}(-m) \surj E_{T_i},
$$
giving a map $T_i\to \SH$. And using the quotient $q_{T_i}$ and
isomorphism $\beta_{T_i}$ we have another family of quotients
parametrized by $T_i$
$$
(V^{\otimes s})^{\oplus c}\otimes \Big(\pi^{}_{T_i*}\big(
\overline{\Delta_i}^* \SP^{\otimes b} \otimes 
\overline{u_{T_i}}^* \SD
\otimes \pi^*_X\SO^{}_X(sm) \big)\Big)^\vee \surj N
$$
Then, using the representability properties of $\SH$ and $P$,
we obtain a morphism to $\SH\times P$, and by lemma \ref{zeroes}
this morphism factors through $Z'$ and since  
a $\delta$-semistable \tensor\ gives a GIT-semistable point (theorem
\ref{thm3.2}), the
image is in $Z^{ss}$. Composing with the 
geometric quotient to $\FM_\delta$ we obtain maps
$$
\hat{f}_i:T_i \stackrel{f_i}\too Z^{ss} \too \FM_\delta 
$$
The morphism $f_i$ is independent of the choice of 
isomorphism $\beta_{T_i}$.
A different choice of isomorphism $g_{T_i}$ will change
$f_i$ to $h_i\cdot f_i$, where $h_i:T_i \to \glv$, so
$\hat{f}_i$ is independent of the choice of $g_{T_i}$.
Then the morphisms $\hat{f}_i$ glue to give a morphism
$$
\hat{f}:T\too \FM_\delta,
$$
and hence we have a natural transformation 
$$
{\SM_\delta} \to \underline{\FM_\delta}.
$$
Recall there is
a tautological family (\ref{3.5bis}) 
of \tensors\ parametrized by $Z'$.
By restriction to $Z^{ss}$, we obtain a tautological family
of  $\delta$-semistable
\tensors\ parametrized by $Z^{ss}$.
If ${\SM_\delta} \to \underline{Y}$ is another natural
transformation, this tautological family defines a 
map $Z^{ss}\to Y$, this factors through the quotient 
$\FM_\delta$, and it is easy to see that this proves that 
$\FM_\delta$ corepresents the functor ${\SM_\delta}$.

Note that in \cite{H-L2}, the moduli space of 
stable framed modules is a fine moduli space. In our situation
this is not true in general, because the analog of the 
uniqueness result of 
\cite[lemma 1.6]{H-L2} does not hold in general for \tensors.
\end{proof}

Now we will give a criterion for S-equivalence. This is 
very similar to the criterion given in \cite{G-S} for
conic bundles.
If $(E,\varphi,u)$ and $(F,\psi,u)$ are two $\delta$-stable
\tensors\ then we have seen that they correspond to the same point
in the moduli space if and only if they are isomorphic. But
if they are strictly $\delta$-semistable, it could happen that
they are S-equivalent (i.e. they correspond to the same point
in the moduli space), even if they are not isomorphic.
Given a \tensor\ $(E,\varphi,u)$, we will construct a canonical
representative of its equivalence class $(E^S,\varphi^S,u)$,
hence $(E,\varphi,u)$ will be S-equivalent to $(F,\psi,u)$ 
if and only if $(E^S,\varphi^S,u)$ is isomorphic to $(F^S,\psi^S,u)$.

Let $(E,\varphi,u)$ be strictly $\delta$-semistable, and let 
$(E_\bullet,m_\bullet)$ be an \textit{admissible} weighted filtration,
i.e. such that
$$
\Big(\sum_{i=1}^t m_i\big( r P_{E_i} -r_i P \big)\Big) +
\mu(\varphi,E_\bullet,m_\bullet) \, \delta \;=\; 0
$$
Let $\SI_0$ be the set of pairs $(k,I)$ 
where $1\leq k \leq c$ is an integer, and
$I=(i_1,\dots,i_s)$ is a multi-index with $1\leq i_j\leq t+1$,
such that the restriction of $\varphi$
$$
\varphi_{k,I}: \overbrace{0 \oplus \ldots  \oplus 0}^{k-1}  
\oplus \big(E_{i_1}\otimes \ldots \otimes E_{i_s}\big) 
\oplus \overbrace{0 \oplus \ldots \oplus 0}^{c-k}
\too (\det E)^{\otimes b} \otimes \fixed_u
$$
is nonzero, and 
$$
\gamma_{r_{i_1}}+\dots+\gamma_{r_{i_s}}=\mu(\varphi,E_\bullet,m_\bullet).
$$
If $(k,I)\in \SI_0$
and $I'=(i'_1,\ldots,i'_s)$ is a multi-index with $I'\neq I$ and 
$i'_j\leq i^{}_j$ for all $j$, then
\begin{equation}
\label{kernel}
\varphi_{k,I'} =0,
\end{equation}
by definition of $\mu(\varphi,E_\bullet,m_\bullet)$.
Hence, if $(k,I)\in \SI_0$, the restriction $\varphi_{k,I}$
defines a homomorphism in the quotient
$$
\varphi'_{k,I}:
 \overbrace{0 \oplus \ldots  \oplus 0}^{k-1}  
\oplus \big(E'_{i_1}\otimes \ldots \otimes E'_{i_s}\big) 
\oplus \overbrace{0 \oplus \ldots \oplus 0}^{c-k}
\too (\det E)^{\otimes b} \otimes \fixed_u,
$$
where $E'_i=E_i/E_{i+1}$.
If $(k,I)\neq \SI_0$, then define $\varphi'_{k,I}=0$.
Finally, we define
$$
E' = E'_1\oplus \ldots \oplus E'_{t+1} \; , \quad 
\varphi'=\bigoplus_{(k,I)} \varphi'_{k,I} \; .
$$
In the definition of $\varphi'$ we are using the 
fact that $\det E\isom \det E'$,
hence $(E',\varphi',u)$ is well-defined up to isomorphism,
and it is called the \textit{admissible deformation} associated
to the admissible filtration $(E_\bullet,m_\bullet)$ of $E$.
Note that it depends on the admissible weighted filtration chosen. 

\begin{proposition}
\label{sequiv}
The \tensor\ $(E',\varphi',u)$ is strictly $\delta$-semistable,
and it is S-equivalent to $(E,\varphi,u)$. If we repeat this
process, after a finite number of iterations the process will
stop, i.e. we will obtain \tensors\ isomorphic to each other.
We call this \tensor\ $(E^S,\varphi^S,u)$
\begin{enumerate}
\item The isomorphism class of $(E^S,\varphi^S,u)$ is independent
of the choices made, i.e. the weighted filtrations chosen.
\item Two \tensors\ $(E,\varphi,u)$ and $(F,\psi,u)$ 
are S-equivalent if and only if $(E^S,\varphi^S,u)$ is
isomorphic to $(F^S,\psi^S,u)$.
\end{enumerate}
\end{proposition}

\begin{proof}
We start with a general observation about GIT quotients. Let $Z$ be a
projective variety with an action of a group $G$ linearized on 
an ample line bundle $\SO_Z(1)$.
Two points in the open subset $Z^{ss}$ of semistable points are
GIT-equivalent (they are mapped to the same point in the moduli space)
if there is a common closed orbit in the closures (in $Z^{ss}$) of
their orbits. Let $z\in Z^{ss}$. Let $B(z)$ be the unique closed orbit in
the closure $\overline{G\cdot z}$ in $Z^{ss}$ of its orbit $G \cdot
z$. Assume that $z$
is not in $B(z)$.
There exists a one-parameter 
subgroup $\lambda$ such that the limit $z_0=\lim_{t\to 0}
\lambda (t) \cdot z$ is in $\overline{G
\cdot z}\setminus G \cdot z$ (for instance, we can 
take the one-parameter subgroup given 
by \cite[Lemma 1.25]{Si}) . Note that we must have
$\mu(z,\lambda)=0$ (otherwise $z_0$ would be unstable). 
Conversely, if $\lambda$ is a one-parameter subgroup
with $\mu(z,\lambda)=0$, then the limit is GIT-semistable
(\cite[Prop. 2.14]{G-S}).
Note that $G \cdot z_0 \subset \overline{G
\cdot z}\setminus G \cdot z$, and then $\dim G \cdot z_0 < \dim G
\cdot z$. Repeating this process with $z_0$ we then get a sequence 
of points that eventually stops and gives $\tilde z \in B(z)$. Two points
$z_1$ and $z_2$ will then be GIT-equivalent if and only if 
$B(z_1)=B(z_2)$.

Let $(E,\varphi,u)$ be a $\delta$-semistable tensor with
an isomorphism $f:V\isom H^0(E(m))$, and let $z=(q,[\Phi])\in Z$
be the corresponding GIT-semistable point. Recall from
proposition \ref{bijection} that there is a bijection between
one-parameter subgroups of $\slv$ 
with $\mu(z,\lambda)=0$ on the one hand,
and weighted filtrations $(E_\bullet,m_\bullet)$
of $E$ with 
$$
\Big(\sum_{i=1}^t m_i(r P_{E_{i}}-r_i P) \Big)+
\mu(\varphi,E_\bullet,m_\bullet)\, \delta \;=\; 0
$$
together with a splitting of the filtration $H^0(E_\bullet(m))$
of $V=H^0(E(m))$ on the other hand.

The action
of $\lambda$ on the point $z$ defines a morphism
$\CC^*\to R_3$ that extends to
$$
h:T=\CC \too Z \, ,
$$
with $h(t)=\lambda(t)\cdot z$ for $t\neq 0$
and $h(0)=\lim_{t\to 0} \lambda(t)\cdot z=z_0$.  

Pulling back the universal family parametrized by
$Z$ by $h$ we obtain the family $(q_T,E_T,\varphi_T,u)$
$$
E^{}_T=\bigoplus_n E^{}_n\otimes t^n \subset 
E\otimes_\CC t^{-N}\CC[t] \subset
E\otimes_\CC \CC[t,t^{-1}]
$$
$$
\begin{array}{rcccc}
q^{}_T: V\otimes\SO^{}_X(-m)\otimes \CC[t] & \stackrel{\gamma}\too &
\oplus_n V^{}_n\otimes\SO^{}_X(-m)\otimes t^n & \too  & E^{}_T \\
 {v^n\otimes 1} & \longmapsto &
{v^n\otimes t^n} & \longmapsto  & {q(v^n)\otimes t^n}
\end{array}
$$
\begin{eqnarray*}
\varphi^{}_T: (E^{}_T{}^{\otimes s})^{\oplus c} 
& \too & 
(\det E_T)^{\otimes b} \otimes 
\overline{u_T}^* \sfixed \otimes \pi_T^*N \\
(\underbrace{0 ,\ldots,0}_{k-1}, 
w_{i_1} t^{i_1}  \cdots w_{i_s} t^{i_s}, 
\underbrace{0 , \ldots , 0}_{c-k})
& \longmapsto &
\varphi(
\underbrace{0 , \ldots  , 0}_{k-1}  
, w_{i_1} \cdots  w_{i_s}, 
 \underbrace{0 , \ldots , 0}_{c-k}
)  \otimes t^{i_1+\cdots+ i_s}
\end{eqnarray*}
Then, as in \cite[\S 4.4]{H-L2}, $(q_t,E_t,\varphi_t,u)$ 
corresponds to $h(t)$ (in particular, if $t\neq 0$, then
$(E_t,\varphi_t,u)$ is canonically 
isomorphic to $(E,\varphi,u)$), and $(E_0,\varphi_0,u)$ is
the admissible deformation associated to 
$(E_{\bullet},m_\bullet)$.

\end{proof}

\section{Orthogonal and symplectic sheaves}
\label{secorthsheaves}

In this section we apply the general theory of \tensors\ to
construct the moduli space of semistable orthogonal 
and symplectic sheaves. The only difference between these
is whether the bilinear form is symmetric or skewsymmetric,
hence we will first consider the orthogonal case, and
at the end of the section we will add some comments
about the symplectic case.
We fix $\fixed_u$ to be $\SO_X$ (i.e. $\param$ is one point and
$\sfixed$ is $\SO_{X\times u}$, and hence we can drop $u$ from 
the notation of tensors).

\begin{definition}
\label{defpgs}
An orthogonal sheaf is a \tensor\ 
$$
(E,\varphi)\; , \qquad \varphi:E\otimes E \too \SO_X
$$
such that
\begin{itemize}

\item (OS1) $(\det E)^{\otimes 2} \isom \SO_X$

\item (OS2) $\varphi$ is symmetric

\item (OS3) $E$ is torsion free

\item (OS4) $\varphi$ induces an isomorphism 
$E|^{}_U\to E|^\vee_U$ 
on the open subset $U$ where $E$ is locally free.
\end{itemize}
An isomorphism of orthogonal sheaves is an isomorphism as tensors.
\end{definition}

It is easy to see that, assuming (OS1) and (OS3), 
the last condition is equivalent to

\begin{itemize}
\item \textit{(OS4$'$)} The induced homomorphism 
$\det E \too \det E^\vee $ is nonzero
(hence an isomorphism).
\end{itemize}
The following lemma justifies this definition for orthogonal sheaves.

\begin{lemma}
\label{bijection1}
There is a bijection between the
set of isomorphism classes of orthogonal
sheaves with $E$ locally
free and the set 
of isomorphism classes of principal $\orth$-bundles.
\end{lemma}

\begin{proof}
The category of principal $\orth$-bundles is equivalent to the category
whose objects are pairs $(P,\sigma)$ (where $\pi:P\to X$ 
is a principal $\glr$-bundle, $\sigma$ is a section of the 
associated fiber bundle
$P\times_{\glr} \glr/\orth$) and whose isomorphisms are isomorphisms
$f:P\to P'$ of principal bundles respecting $\sigma$ (i.e.
$\pi'\circ f=\pi$ and $(f\times \id_{\orth}) \circ \sigma = \sigma'$).
Note that this notion of isomorphism is \textit{not} the same
as isomorphism of reductions.

The category of principal $\glr$-bundles is equivalent to the
category of vector bundles of rank $r$.
The quotient $\glr/\orth$ is the set of invertible
symmetric matrices (send $A\in \glr$ to $({}^T\! A^{-1} A^{-1})$).
Hence, a section $\sigma$ is the same 
thing as a homomorphism $\varphi$ as in (OS2).
Now it is easy to check that there is a bijection betweeen these 
sets of isomorphisms classes.

\end{proof}

\begin{remark}\textup{ 
Note that the categories are not
equivalent: for example, let $P$ be a simple principal $G$-bundle,
i.e. the set of automorphisms of $P$ is the center of $G$ (a finite
group), but the set of automorphisms of the corresponding 
$G$-sheaf is $\CC^*$. We will have an equivalence 
of categories if we consider only
isomorphisms $(f,\alpha)$ with $\alpha=1$, as in remark \ref{alpha1}. 
This would be important
if we wanted to construct the moduli stack, but since we are
interested in the moduli space this is irrelevant, because 
the moduli space does not detect the group of automorphisms.
}
\end{remark}

Let $(E,\varphi)$ be an orthogonal (or symplectic) sheaf.
A subsheaf of $F$ of $E$ is called isotropic 
if $\varphi|_{F\otimes F}=0$.
Given a subsheaf $i:F\inj E$, using the bilinear form $\varphi$
we can associate
the perpendicular subsheaf
$$
F^\perp =\ker (E \stackrel{\overline{\varphi}}{\too} 
E^\vee \stackrel{\overline{i^\vee}}{\too} F^\vee),
$$
where $\overline{\varphi}:E\to E^\vee$ is the homomorphism induced by
$\varphi$.

\begin{definition}[Stability]
\label{principalst}
An orthogonal sheaf $(E,\varphi)$ is (semi)stable 
if for all isotropic subsheaves $F\subset E$, 
\begin{equation}
\label{eqcritGM}
 P^{}_{F} + P_{F^\perp}   \;(\preceq)\; P.
\end{equation}
An orthogonal sheaf $(E,\varphi)$ is 
slope-(semi)stable if 
for all isotropic
subsheaves $F\subset E$, 
$$
\deg(F)\;(\leq)\; 0.
$$
\end{definition}

As usual, we can assume that $F$ is saturated.
A family of semistable orthogonal sheaves 
parametrized by $T$ is a family of \tensors\ 
\begin{equation}
\label{familyorth}
(E_T,\varphi_T,N), \quad 
\varphi_T:E_T\otimes E_T\too \pi_T^* N,
\end{equation}
such that $(\det E_T)^{\otimes 2}$ is 
isomorphic to the pullback of some line bundle on $T$, $\varphi_T$
is symmetric, and $\varphi_T$ induces
an isomorphism 
$E_T|_U \to E_T^\vee\otimes \pi_T^* N|_U$ on the open
set $U$ where $E_T$ is locally free, and such that the restriction
to $X\times t$ for all closed points $t$ is a semistable
orthogonal sheaf. 

Using this notion of family, we
define the functor $\SM_{\orth}$   
of semistable orthogonal sheaves. 
We will construct 
a moduli space corepresenting this functor (theorem \ref{projective}).

In proposition \ref{critGM} we show that an
orthogonal (or symplectic) sheaf 
$(E,\varphi)$ is (semi)stable in this sense if and only if
it is $\delta$-(semi)stable as a \tensor\ (definition
\ref{stability}), provided that 
$\delta_1>0$. Hence, the moduli space of semistable 
orthogonal (or symplectic) sheaves 
is a subscheme of the moduli space of $\delta$-semistable \tensors.
In theorem \ref{projective} we show that it is in fact projective.
We can also ask about slope-semistability, and in proposition 
\ref{critslope} we show that slope-(semi)stability in this sense and
slope-$\tau$-(semi)stability as a \tensor\ coincide if $\tau>0$.
If $\delta_1=0$, then the notion of $\delta$-semistability 
as a \tensor\ is not
equivalent to semistability as an orthogonal sheaf. At the end of 
the section we give an example of this.

We start with some preliminaries.
The intersection $F\cap F^\perp$ is an isotropic subsheaf of $F$.
The following lemma gives exact sequences relating these subsheaves.

\begin{lemma}
\label{short}
With the previous notation:
\begin{enumerate}

\item Let $U$ be the open set where $F$, $E$ and $E/F$ are locally
free. There is an exact sequence on $U$
$$
0 \too F^\perp|_U \too E|_U \too 
F^\vee|_U \too 0,
$$
and hence $\rk(F^\perp)=\rk(E)-\rk(F)$. If furthermore $F$ is
saturated (i.e. $E/F$ is torsion free), then $\codim(X-U)\geq 2$
and hence $\deg(F^\perp)=\deg(F)$.

\item If $F$ is saturated, then $F\cap F^\perp$ is also saturated.

\item There is an exact sequence
\begin{equation}
\label{cap}
0 \too F\cap F^\perp \too F\oplus F^\perp \too F+F^\perp \too 0
\end{equation}
\item $F+F^\perp \subset (F\cap F^\perp)^\perp$, 
$\rk(F+F^\perp) = \rk((F\cap F^\perp)^\perp)$, and hence
$\deg(F+F^\perp) \leq \deg((F\cap F^\perp)^\perp)$.

\item Let $F$ be a saturated subsheaf. If $F\cap F^\perp \neq 0$, then
\begin{equation}
\label{shortA}
\deg(F)\leq \deg(F\cap F^\perp),
\end{equation}
and if $F\cap F^\perp=0$, then 
\begin{equation}
\label{shortB}
\deg(F)\leq 0.
\end{equation}

\end{enumerate}
\end{lemma}

\begin{proof}
Since $E/F|_U$ is locally free, the last term in the following exact
sequence is zero
$$
0 \too (E/F)^\vee|_U \too E^\vee|_U 
\stackrel{i^\vee|_U}{\too} F^\vee|_U \too
Ext^1((E/F)|_U,\SO_U)=0,
$$
and hence $i^\vee|_U$ is surjective. 
Combining this with (OS4) we get the exact sequence
$$
0 \too F^\perp|_U \too E|_U \cong E^\vee|_U
\stackrel{i^\vee|_U}{\too} 
F^\vee|_U \too 0.
$$
If $E/F$ is torsion free, then
$\codim(X-U)\geq 2$ and we can use this sequence 
to obtain $\deg(F^\perp)=\deg(F)$.

To prove item 2, first we show that $F^\perp$ is saturated.
The composition $E\to E^\vee\to F^\vee$ factors as
$$
E \surj E/F^\perp \inj F^\vee.
$$
The sheaf $F^\vee$ is torsion free, 
and hence also $E/F^\perp$ is torsion free.

We conclude by showing that the stalk $(E/(F\cap F^\perp))_x
=E_x/(F_x\cap F_x^\perp)$ is torsion free for all points $x\in X$.
Let $v \in E_x$ and let $0\neq f\in \mathfrak{m}_x$ be a nonzero element
in the maximal ideal of the local ring of $x$, 
such that $fv \in F^{}_x\cap F_x^\perp$.
Since $fv \in F_x$, and $F_x$ is saturated, 
then $v\in F_x$. The same argument applies to $F_x^\perp$, and hence
$v\in F_x\cap F_x^\perp$.

Items 3 and 4 are easy to check. To show item 5, if 
$F\cap F^\perp\neq 0$, use the exact sequence (\ref{cap}),
together with items 1, 2 and 4. If $F\cap F^\perp=0$, then
$F\oplus F^\perp = F+F^\perp$ is a subsheaf of $E$ of rank $r$, then 
$\deg(F)+\deg(F^\perp)\leq 0$, and hence $\deg(F)\leq 0$.

\end{proof}

The fact that on a generic fiber the quadratic form is nondegenerate
has the following useful consequence:

\begin{lemma}
\label{negative}
If $(E,\varphi)$ is an orthogonal or syplectic sheaf, then
for all weighted filtrations
\begin{equation}
\label{mu0}
\mu(\varphi,E_\bullet,m_\bullet)\leq 0.
\end{equation}
\end{lemma}

\begin{proof}
First we will show that if $Q:W\otimes W \to \CC$ 
is a bilinear
nondegenerate form on a vector space $W$, then $\overline{Q}
\in \PP(W^\vee\otimes W^\vee)$ is GIT-semistable under
the natural action of $\slw$ (with the natural linearization 
induced on $\SO(1)$). 
The point $\overline{Q}$ is unstable if and only if there is
a one-parameter subgroup $\lambda$ of $\slw$ such that 
$\lim_{t\to 0} \lambda(t)\cdot Q=0$. But this is impossible because 
$\det(\lambda(t))=1$, and then
$$
\det(\lambda(t)\cdot Q)= \det(\lambda(t)\, Q \;{}^T\!\lambda(t)) = 
\det(Q) \neq 0,
$$
hence $\overline{Q}$ is semistable.
Then, using this and condition (OS4), it follows that 
$$
\mu(\varphi,E_\bullet,m_\bullet)\leq 0
$$
for all weighted filtrations. 
\end{proof}

\begin{proposition}
\label{critGM}
Assume $\delta_1>0$. An orthogonal 
sheaf $(E,\varphi)$ 
is (semi)stable
if and only if it is $\delta$-(semi)stable as a \tensor.
\end{proposition}

\begin{proof}
To see that $\delta$-(semi)stable as a \tensor\ implies (semi)stable as 
an orthogonal
sheaf, we apply the stability
condition to the weighted
filtration $F\subset F^\perp \subset E$ with weights $m_1=m_2=1$.
By lemma \ref{short}(1), $r=\rk(F)+\rk(F^\perp)$.
Since $F$ is isotropic, $\mu(\varphi,E_\bullet,m_\bullet)=0$,
hence the stability condition (\ref{stformula}) gives the result:
$$
r\big( P^{}_{F} + P_{F^\perp}  - P \big) =
\big(r  P^{}_{F} - \rk(F)P\big) + 
\big(r P_{F^\perp}  -\rk(F^\perp)P\big) 
 \;(\preceq)\; 0.
$$

Now we will show that if $(E,\varphi)$ is (semi)stable as an
orthogonal sheaf, then it is $\delta$-(semi)stable as a \tensor.
We start with a vector space $W$ and a nondegenerate bilinear
form $Q:W\otimes W\to \CC$. Let $(W_\bullet, m_\bullet)$ be a 
weighted filtration with
\begin{equation}
\label{muQ}
\mu(Q,W_\bullet, m_\bullet)=0.
\end{equation}
Denote $r_i=\dim W_i$.
Take a basis of $W$ adapted to the filtration, and let $\lambda$
be the one-parameter subgroup of $\slw$ associated to this basis
and weights $m_\bullet$. Let $\gamma=\sum_{i=1}^{t} m_i\gamma^{(r_i)}$
as in (\ref{gammak}).
Since $\mu(Q,W_\bullet, m_\bullet)=0$, the limit $Q'= 
\lim_{t\to 0}\lambda(t)\cdot Q$
exists, and $\det Q' =\det Q$.
Furthermore, we also have 
\begin{equation}
\label{muQ0}
\mu(Q',W_\bullet, m_\bullet)=0.
\end{equation}
Write $Q$ and $Q'$ as block matrices
$$
Q= \left(
\begin{array}{cccc}
Q_{1,1} & Q_{1,2} & \dots & Q_{1,t+1} \\
Q_{2,1} & Q_{2,2} & \dots & Q_{2,t+1} \\
\ldots & \ldots & \dots & \ldots \\
Q_{t+1,1} & Q_{t+1,2} & \dots & Q_{t+1,t+1} \\
\end{array}
\right)
$$
Note that if $\gamma_{r_i}+\gamma_{r_j}<0$, then $Q_{i,j}=0$
because of (\ref{muQ}).
We have
\begin{equation}
\label{weights}
Q'_{i,j}= \left\{
\begin{array}{cl}
0       &, \,\gamma_{r_i}+\gamma_{r_j}<0 \\  
Q_{i,j} &, \,\gamma_{r_i}+\gamma_{r_j}=0 \\
0       &, \,\gamma_{r_i}+\gamma_{r_j}>0 \\
\end{array}
\right .
\end{equation}
The weights $\gamma_{r_i}+\gamma_{r_j}$ strictly increase with both $i$ 
and $j$. Assume $Q'_{i,j}\neq 0$. Then, if $(a,b)\neq(i,j)$,
and either $a\leq i$, $b\leq j$, or $a\geq i$, $b\geq j$, we have    
$Q'_{a,b}=0$. In matrix form:
\begin{equation}
\label{zeros}
Q'= \left(
\begin{array}{cccc}
0       &    0    &       &      \\
0       & Q'_{i,j} &  0    & 0    \\
        &    0    &  0    & 0    \\
        &    0    &  0    & 0    \\
\end{array}
\right)
\end{equation}
Since $\det Q'=\det Q\neq 0$,
in each row of $Q'$ there 
must be at least one nonzero block (and the same for columns).
This, together with (\ref{zeros}) implies
\begin{equation}
\label{q0}
Q'= \left(
\begin{array}{ccccc}
   0      &    0    & \ldots &   0      & Q'_{1,t+1}  \\
   0      &    0    & \ldots & Q'_{2,t}  &     0      \\
\ldots    & \ldots  & \ldots & \ldots   &    \ldots  \\
   0      & Q'_{t,2} & \ldots &   0      &     0      \\
Q'_{t+1,1} &    0    & \ldots &   0      &     0      \\
\end{array}
\right)
\end{equation}
with nonzero blocks in the second diagonal, and zero everywhere
else. 
Since $Q'$ is nondegenerate, these blocks give isomorphisms
for all $1\leq i\leq t+1$
$$
Q'_{i,t+2-i}: W_{t+2-i}/W_{t+1-i} \stackrel{\cong}{\too} 
 W_i/W_{i-1},
$$
and a short calculation then gives $r_{i}=r-r_{t+1-i}$.
This, together with (\ref{q0}), implies that 
\begin{equation}
\label{perpw}
W_i^\perp=W^{}_{t+1-i}.
\end{equation}
Finally (\ref{weights}) and (\ref{q0}) imply that 
$\gamma_{r_i}+\gamma_{r_{t+2-i}}=0$ for all $1\leq i \leq t+1$.
Then, using this and the definition of $\gamma$,
$$
0 = 
(\gamma_{r_{i+1}}+\gamma_{r_{t+1-i}})- 
(\gamma_{r_i}+\gamma_{r_{t+2-i}})=
r(m_i-m_{t+1-i}).
$$

Let $(E_\bullet,m_\bullet)$ be a weighted filtration with
$\mu(\varphi,E_\bullet,m_\bullet)=0$.
We can assume that all subsheaves $E_i$ are saturated.
Apply the previous argument to $W=E|_x$, the fiber
over a point where $E$ is locally free, and $Q$
the bilinear form induced by $\varphi$ on the
fiber.
We have $(\ref{perpw})$, hence it follows that 
$E_i^\perp \supset E^{}_{t+1-i}$.
Furthermore, as we have just seen 
$m_i=m_{t+1-i}$ and $r_{i}=r-r_{t+1-i}$ for all $i$. 
Hence we can write
\begin{eqnarray*}
\sum_{i=1}^t m_i\big( r P_{E_i} -r_i P \big) \;=\;   \\
\sum_{i=1}^{[(t+1)/2]} m_i r \big(  P_{E_i} + P_{E_{t+1-i}} - P \big) 
\;\preceq\;  \\
\sum_{i=1}^{[(t+1)/2]} m_i r \big(  P^{}_{E_i^{}} + P_{E_i^\perp} - P \big) 
\;(\preceq)\; 0,
\end{eqnarray*}
where the last inequality is given by (\ref{eqcritGM}).

Let $(E_\bullet,m_\bullet)$ be a weighted filtration with
$\mu(\varphi,E_\bullet,m_\bullet)\neq 0$ . By
lemma \ref{negative}, it is strictly negative.

We claim that $\deg(E_i)\leq 0$ for all $i$. Assume that this
is not true. Then there is a saturated subsheaf $F\subset E$ with
$\deg(F)>0$. By lemma \ref{short}(5), $N=F\cap F^\perp\neq 0$, 
and $0<\deg(F)\leq \deg(N)$. By lemma \ref{short}(2), $N$ is
saturated, and by lemma \ref{short}(1), $\deg(N)=\deg(N^\perp)$.
Consider the weighted filtration $N\subset N^\perp \subset E$
with weights $m_1=m_2=1$. Since $N$ is isotropic, 
$\mu(\varphi,N\subset N^\perp)=0$, and since
$\deg(N)=\deg(N^\perp)>0$, this weighted filtration 
contradicts (\ref{eqcritGM}).

Hence, using $\deg(E_i)\leq 0$ together with $\delta_1>0$,
\begin{eqnarray*}
\Big(\sum_{i=1}^t m_i\big( r P_{E_i} -r_i P \big)\Big) +
\mu(\varphi,E_\bullet,m_\bullet) \, \delta \;=   \hspace{2cm} \\
\Big(\big(\sum_{i=1}^t m_i\frac{ r \deg(E_i)}{(\di-1)!}  \big) +
\mu(\varphi,E_\bullet,m_\bullet) \, \delta_1 \Big) t^{\di-1} +
O(t^{\di-2}) \;\prec \; 0
\end{eqnarray*}

\end{proof}

\begin{proposition}
\label{critslope}
Assume $\tau>0$.
An orthogonal
sheaf $(E,\varphi)$ 
is slope-(semi)stable
if and only if it is slope-$\tau$-(semi)stable as a \tensor.
\end{proposition}

\begin{proof}
The proof of proposition \ref{critGM}, replacing the 
Hilbert 
polynomials
$P_{F}$, $P_{E_i}$, $P$,...  by the degrees 
$\deg(F)$, $\deg(E_i)$, $d$,... proves that 
$(E,\varphi)$ is slope-$\tau$-(semi)stable if and only if
for all isotropic subsheaves $F\subset E$,
$$
\deg(F) + \deg(F^\perp) \leq \deg(E).
$$
We can assume that $F$ is saturated, hence $\deg(F)
=\deg(F^\perp)$ by lemma \ref{short}(1), and since $\deg(E)=0$, 
the result follows.
\end{proof}

Fix a polynomial $P$. Recall that $\SM_{\orth}$ is the 
functor of families of semistable orthogonal sheaves.
Define $\FM_{\orth}$ 
to be the subscheme of the moduli space of $\delta$-semistable 
\tensors\ corresponding
to orthogonal
sheaves with Hilbert polynomial $P$.
The notion of S-equivalence for orthogonal sheaves is the same that
was described in proposition \ref{sequiv}.

\begin{theorem}
\label{projective}
The scheme $\FM_{\orth}$ is a coarse moduli space of
S-equivalence classes of semistable orthogonal sheaves.
There is an open subscheme $\FM^0_{\orth}$ corresponding to semistable 
orthogonal bundles.
\end{theorem}

\begin{proof}
The proof that $\FM_{\orth}$ corepresents the functor
$\SM_{\orth}$ is completely analogous to the
proof of theorem \ref{mainthm} (see section \ref{sectheorem}),
so we will not repeat it.
The subscheme $\FM^0_{\orth}$ is open because being locally free is an open
condition.

Now we will prove that this moduli space is projective.
Conditions (OS1) and (OS2) are closed conditions, so they 
define a projective subscheme $\FM_{1,2}$ of the moduli space of 
$\delta$-semistable \tensors.
The lemma will be proved by showing that 
$\FM_{\orth}=\FM_{1,2}$.
If $(E,\varphi)$ is $\delta$-semistable then $E$ is torsion free,
so it only remains to check that if condition (OS4) does not
hold, then $(E,\varphi)$ is $\delta$-unstable.

Assume that the homomorphism $\det E \to \det E^\vee $ induced by
$\varphi$ is zero.
Then the sheaf $E^\perp$ defined as
$$
0 \too E^\perp \too E \stackrel{\overline{\varphi}}{\too} E^\vee
$$ 
is nonzero. Let $C$ be the cokernel of $\overline{\varphi}$
$$
E \stackrel{\overline{\varphi}}{\too} E^\vee \too C \too 0.
$$
Taking the dual of this sequence and restricting to the open subset $U$ 
of $X$ where $E$ is locally free, we get
$$
0 \too C^\vee|_U \too E^{\vee\vee}|_U=E|_U 
\stackrel{{\overline{\varphi}}|^\vee_U}{\too} E^\vee|_U
$$
By (OS2) we have 
$\overline{\varphi}|^{}_U= \overline{\varphi}|^\vee_U$,
hence 
$\ker(\overline{\varphi}|^{}_U) \cong \ker(\overline{\varphi}|^\vee_U)$,
and then $C^\vee|_U \cong E^\perp|_U$, and since $\codim(X-U)\geq 2$, 
$\deg(C)=-\deg(E^\perp)$. The exact sequence on $U$
$$
0 \too E^\perp|_U \too E|_U \too E^\vee|_U 
\too C|_U \too 0
$$
implies that $\deg(E^\perp)=0$.
Consider the weighted filtration
$0\subset E^\perp \subset E$, $m_1=1$. We have
$$
\mu(\varphi,E_\bullet,m_\bullet)>0.
$$
Recall that $\tau=\delta_1 (\di-1)!$. Then 
$$
r \deg(E^\perp)-\rk(E^\perp)\deg(E) +
\mu(\varphi,E_\bullet,m_\bullet) \tau =
\mu(\varphi,E_\bullet,m_\bullet) \tau >0,
$$
and hence $(E,\varphi)$ is slope-$\tau$-unstable
(definition \ref{slopestability}), and in particular,
$\delta$-unstable.

\end{proof}

\begin{remark}
\textup{
The same proof gives that if $(E,\varphi)$ is a slope-$\tau$-semistable
\tensor\  with $\tau>0$, satisfying conditions (OS1),
(OS2) and (OS3), then condition (OS4) holds.
}
\end{remark}

\medskip
\noindent\textbf{Example.}
We will give an example showing that, if we do not
require $\delta_1$ to be positive, 
the notion of $\delta$-stability as a \tensor\ (definition \ref{stability})
is different from the notion of stability as an orthogonal sheaf
(definition \ref{principalst}). We will check this by showing 
an example of an orthogonal sheaf whose $\delta$-stability really
depends on  $\delta$.

Let $X=\PP^2$, let $p^{}_1$, $p^{}_2$, $p^{}_3$ be three different
points in $\PP^2$, and consider the ideal sheaves
$I_1=I_{p^{}_1}$ and $I_2=I_{p^{}_2\cup p^{}_3}$. Let 
$$
(E,\varphi)= 
\left(I_2 \oplus I_1 \oplus \SO_X, 
\left(
\begin{array}{ccc}
0 & 0 & 1/2 \\
0 & 1 & 1/2 \\
1/2 & 1/2 & 1 \\
\end{array}
\right)
\right)
$$
In particular, the first summand $I_2$ of $E$ is isotropic, and
$I_2^\perp=I^{}_2 \oplus I^{}_1$.
Let $\delta=\delta_1 t+\delta_2 \in \QQ[t]$ be a polynomial 
as in (\ref{delta}).

\begin{lemma}
If $\delta_1=0$ and $0<\delta_2<3/2$, then $(E,\varphi)$ is
$\delta$-unstable as a \tensor. If $\delta_1>0$, then $(E,\varphi)$ is 
$\delta$-stable as a \tensor.
\end{lemma}

\begin{proof}
The first claim is proved by considering the filtration
$\SO_X \subset I_2\oplus I_1 \oplus \SO_X$. If $\delta_1=0$,
then this filtration does not satisfy (\ref{stformula}), 
hence contradicts semistability.

Now we will prove that if $\delta_1>0$, then $(E,\varphi)$ is 
$\delta$-stable. Using proposition \ref{critGM}, 
we only have to study filtrations of the form
$$
E_1 \subset E_2=E_1^\perp \subset E,
$$
with $\rk(E_1)=1$, $\rk(E_2)=2$ and $E_1$ isotropic
and saturated. Using the Riemann-Roch formula we have
$$
P_{E_1} + P_{E_2} - P \;=\; -2 c_2(E_1) -2 c_2(E_2) + 2 c_2(E),
$$
so we need to estimate the second Chern classes of $E_1$ and $E_2$.

The sheaf $E^{}_2=E_1^\perp$ is saturated (see the proof of lemma
\ref{short}(2)). Define the torsion free
rank one subsheaf $J$
$$
0 \too E_2 \too I_2\oplus I_1 \oplus \SO_X \stackrel{(a,b,c)}\too J
\too 0,
$$
where $a$, $b$ and $c$ are respectively elements of $\Hom(I_2,J)$, 
$\Hom(I_1,J)$ and $\Hom(\SO_X,J)$. We have $\deg(J)=0$, so $J$ is 
the ideal sheaf of a zero-dimensional subscheme of $\PP^2$.
We distinguish several cases:

\begin{itemize}

\item If $c\neq 0$, then $J=\SO_X$, and $c_2(E)=3$.

\item If $a\neq 0$, $b\neq 0$, then again $J=\SO_X$, and $c_2(E)=3$.

\item If $a=0$, $b\neq 0$ and $c=0$, then $E_2$ does not contain 
a subsheaf $E_1$ with $E_2= E_1^\perp$, hence this cannot
happen.

\item If $a\neq 0$, $b=0$ and $c=0$, then again $E_2$ does not contain 
a subsheaf $E_1$ with $E_2 \subset E_1^\perp$, hence this cannot
happen.
\end{itemize}
So we conclude that $c_2(E_2)=3$.
The sheaf $E_1$ is a rank one subsheaf of $I_2\oplus I_1 \oplus
\SO_X$, hence $c_2(E_1)>0$ unless $E_1$ is the third summand $\SO_X$,
but this is not possible because the third summand is not
isotropic. Putting everything together,
$$
P_{E_1} + P_{E_2} - P \; =\;
-2c_2(E_1)-2c_2(E_2)+2c_2(E) \;\prec\; 0,
$$
hence $(E,\varphi)$ is $\delta$-stable by proposition \ref{critGM}.
\end{proof}

\begin{remark}\textup{
Note that $(E,\varphi)$ is stable as an orthogonal sheaf, but
$E$ is Gieseker-unstable as a sheaf.}

\textup{
On the other hand, an orthogonal sheaf $(E,\varphi)$ is 
slope-semistable if and only if
$E$ is slope-semistable as a sheaf. Indeed, if $F$ is a saturated 
subsheaf of $E$ with
$\deg(F)>0$, then lemma \ref{short}(5) shows that the isotropic
subsheaf $F\cap F^\perp$ is nonzero and has positive degree, hence
$(E,\varphi)$ is slope-unstable.
}
\end{remark}

To obtain symplectic sheaves instead of orthogonal sheaves,
we only need to take $\varphi$ skewsymmetric instead of
symmetric. It follows that $\det E=\SO_X$ (recall that 
for orthogonal sheaves we only had $(\det E)^{\otimes 2}=\SO_X$

There is a bijection between the set of isomorphism classes of
symplectic bundles and principal $\sympl$-bundles.
The proof is the same as with orthogonal bundles, 
noting that the quotient $\glr/\sympl$ is the set of invertible
antisymmetric matrices (send $A\in \glr$ to 
$({}^T\! A^{-1} J A^{-1})$, where $J$ is the matrix representing the
standard symplectic structure of $\CC^r$). 

All the results for orthogonal sheaves hold for symplectic sheaves,
and in particular there is a coarse moduli space of S-equivalence
classes of semistable symplectic sheaves.

\section{Special orthogonal bundles}
\label{secspecialgbundles}

\begin{definition}[Special orthogonal sheaf]
\label{defpss}
A special orthogonal sheaf is a triple
$$
(E,\; \varphi:E\otimes E \too \SO_X ,\; \psi:\det E\too \SO_X) 
$$
such that
\begin{itemize}
\item (SOS1) $\psi$ is an isomorphism.

\item (SOS2) $\varphi$ is symmetric.

\item (SOS3) $E$ is torsion free.

\item (SOS4) $\varphi$ induces an isomorphism 
$E|^{}_U\to E|^\vee_U$ 
on the open subset $U$ where $E$ is locally free.

\item (SOS5) $\det(\varphi)=\psi^2$. More precisely,
let $\varphi':E\to E^\vee$
and $\psi':\SO_X \to \det E^\vee $ be the associated 
homomorphisms. Then we require $\det(\varphi')=\psi\otimes \psi'$.

\end{itemize}
An isomorphism of special orthogonal sheaves is a pair $(f,\lambda)$
such that $f:E\to E'$ is an isomorphism, $\lambda\in \CC^*$
and the following diagrams commute
\begin{equation}
\label{sorisom}
\xymatrix{
{E\otimes E}  \ar[r]^{f^{\otimes 2}} \ar[d]^{\varphi} &  
{E'\otimes E' } \ar[d]^{\varphi'}
 & {\det E} \ar[r]^{\det f} \ar[d]^{\psi} & 
{\det E' } \ar[d]^{\psi'} \\
{\SO_X} \ar[r]^{\lambda^2}  & {\SO_X} 
 & {\SO_X} \ar[r]^{\lambda^r}  & {\SO_X}
}
\end{equation}
\end{definition}
It is easy to see that, assuming (SOS1) and (SOS3), 
condition (SOS4) is equivalent to

\begin{itemize}
\item \textit{(SOS4$'$)} The induced homomorphism 
$\det E\too \det E^\vee$ is nonzero
(hence an isomorphism).
\end{itemize}
Condition (SOS5) is equivalent to 

\begin{itemize}
\item \textit{(SOS5$'$)}.
Let $U$ be the open subset where $E$ is torsion free.
For all $x\in U$,
fix a basis of the fiber $E_x$ of $E$ on $x$,
Using this basis (and the canonical identification $\SO_x\isom \CC$), 
$\varphi$ restricted to $x$ gives a symmetric
matrix $\varphi(x)$, and $\psi$ restricted to $x$ gives a
complex number $\psi(x)$. Then we require $\det(\varphi(x))=\psi(x)^2$.
\end{itemize} 
This definition of special orthogonal sheaf is justified by the
following lemma.
\begin{lemma}
There is a bijection between the
set of isomorphism classes of special orthogonal sheaves with $E$ locally
free and the set 
of isomorphism classes of principal $\sor$-bundles.
\end{lemma}

\begin{proof}
Let $\sor$ act by multiplication on the right on $\glr$, and
consider the quotient $\glr/\sor$. Let $A\in \glr$, and
let $[A]$ be the class in $\glr/\sor$. To this class we
associate the pair $({}^T\! A^{-1} A^{-1}, \det(A^{-1}))$.
This gives a bijection between the set $\glr/\sor$ 
and the set of pairs
$(B,\beta)$, where $B$ is a symmetric invertible matrix and $\beta$
is a nonzero complex number such that 
$$
\det B=\beta^2.
$$

Given a principal $\glr$-bundle $P$ (or equivalently a vector bundle
$E$), a reduction of structure group to $\sor$ is a section $\sigma$
of the associated bundle $P\times_{\glr} \glr/\sor$, and then this 
is equivalent to a pair of homomorphisms $(\varphi,\psi)$ as in 
definition \ref{defpss}.

The rest of the proof is analogous to the proof of lemma \ref{bijection1}. 
\end{proof}

\begin{definition}[Stability]
\label{stsor}
A special orthogonal sheaf $(E,\varphi,\psi)$ 
is (semi)stable if the
associated orthogonal sheaf $(E,\varphi)$ 
is (semi)stable.
\end{definition}

A family of semistable special orthogonal sheaves parametrized by
$T$ is a tuple $(E_T,\varphi_T,\psi_T,N)$ 
such that $(E_T,\varphi_T,N)$ is a family of semistable 
orthogonal sheaves (cf. (\ref{familyorth})), and
$\psi_T:\det E_T \to \pi^*_T L$ is an isomorphism, where
$L$ is a line bundle on $T$.
Two families are isomorphic if there is a pair $(f,\lambda:M\to M')$
where $f:E^{}_T\to E'_T$ is an isomorphism, 
$M^{\otimes 2}\isom N$, ${M'}^{\otimes 2}\isom N'$, 
$\lambda$ is an isomorphism, and the relative versions of the diagrams
(\ref{sorisom}) commute.
In this section (theorem \ref{modsor}) 
we will construct the moduli space of semistable
special orthogonal sheaves (with fixed Hilbert polynomial).

There is a map between isomorphism classes
$$
\begin{array}{ccc}
 \Big\{
 \textup{Special orthogonal sheaves}
 \Big\} 
& 
 \stackrel{f}{\too}
&
 \left\{
 \begin{array}{c}
 \textup{Orthogonal sheaves}\\
 \textup{such that }\det E\isom \SO_X \\
 \end{array}
 \right\}
\\
 (E,\varphi,\psi) & \longmapsto &(E,\varphi) 
\\
\end{array}
$$
This map will induce a morphism between the corresponding
moduli spaces.

\begin{lemma}
\label{twotoone}
Let $(E,\varphi)$ be an orthogonal sheaf such that 
$\det E\isom \SO_X$. 

If $E$ has an automorphism $f$ such that $f\otimes f=\id_{E\otimes E}$
and $\det f =-\id_{\det E }$, then the preimage of $(E,\varphi)$ 
under the map $f$ 
consists of exactly one isomorphism class.

If $E$ does not have such an automorphism, 
then the preimage consists of exactly two distinct 
isomorphism classes, represented
by two special orthogonal sheaves $(E,\varphi,\psi)$ and 
$(E,\varphi,-\psi)$, differing in the sign
of the isomorphism $\psi$.
\end{lemma}

\begin{proof}
Property (SOS5) implies that to obtain the isomorphism $\psi$
we have to extract a square root, so we obtain two 
special orthogonal sheaves $P=(E,\varphi,\psi)$ and $P'=(E,\varphi,-\psi)$
mapping to the given orthogonal sheaf. It only remains to
check if theses two objects are isomorphic or not.

If there is an automorphism $f:E\to E$ with the above properties,
then $(f,1)$ is an isomorphism between $P$ and $P'$.

Conversely, assume that there is an isomorphism $(f,\lambda)$ between
$P$ and $P'$. Then $f'=(1/\lambda)f$ is an automorphism of $E$ with
$f'\otimes f'=\id$ and $\det f'=-\id$.
\end{proof}

\begin{corollary}
If $r$ is odd, there is a bijection between the set of 
isomorphism classes of special orthogonal sheaves and the 
set of isomorphism classes of orthogonal sheaves 
with $\det E \isom \SO_X$. 
\end{corollary}

\begin{proof}
Apply lemma \ref{twotoone} to $f=-\id_E$ (multiplication by $-1$).
\end{proof}

In particular, for $r$ odd, the moduli space of
(semi)stable special orthogonal sheaves consists of the 
components of the moduli space of (semi)stable  
orthogonal sheaves with trivial determinant.
On the other hand, if $r$ is even and $E$ is simple, then
for each orthogonal sheaf with trivial determinant, we
have two nonisomorphic special orthogonal sheaves.
From now on we will assume that $r$ is even.

Fix a Hilbert polynomial $P$.
Let $m$ be a large
integer number as in section \ref{secgitconstruction}. Let $V$ be a
vector space of dimension $P(m)$.
Let $(g,E,\varphi,\psi)$ be a tuple where $(E,\varphi,\psi)$ is
a semistable special orthogonal sheaf and $g$ is an isomorphism
between $H^0(E(m))$ and $V$. As in section \ref{secgitconstruction},
the homomorphism $\varphi$ gives a vector
$$
\Phi \in (V^{\otimes 2})^\vee \otimes H^0(\SO_X(2m))
$$
We denote $\Phi_s=\Phi^{\otimes r/2}$ the associated vector
\begin{equation}
\label{symr2}
\Phi_s \in \Sym^{r/2} \big((V^{\otimes 2})^\vee \otimes H^0(\SO_X(2m))\big).
\end{equation}
The homomorphism $\psi$ induces a linear map
$$
\Psi:\bigwedge{}^r V \too H^0(\det(E)(rm)) \too H^0(\SO_X(rm)),
$$
and hence a vector (that we denote with the same letter) 
$$
\Psi\in (\bigwedge{}^r V)^\vee \otimes H^0(\SO_X(rm)).
$$
These two quotients give a point $[\Phi_s,\Psi]$ in the 
projective space $\wt{P}$ defined as
$$
\PP\Big(
\Sym^{r/2} \big((V^{\otimes 2})^\vee \otimes H^0(\SO_X(2m))\big)
\oplus \big( (\bigwedge{}^r V)^\vee \otimes H^0(\SO_X(rm)) \big) 
\Big).
$$
It is easy to check that the point only depends on the isomorphism
class of the tuple. Here it is crucial that we took the
$r/2$-symmetric power in (\ref{symr2}): take the isomorphism 
$\lambda\id:E\to E$ (multiplication by $\lambda\in \CC^*$). It
sends $\Phi$ to $\lambda^2 \Phi$, and $\Psi$ to $\lambda^r \Psi$,
hence it sends $[\Phi_s,\Psi]$ to $[\lambda^r \Phi_s,\lambda^r
\Psi]$, and this is the same point in the projective space.

Let $\SH$ be the Hilbert scheme of quotients as in section 
\ref{secgitconstruction}, and then given a tuple
$(g,E,\varphi,\psi)$ we associate a point $(q,[\Phi_s,\Psi])$ in
$\SH\times \wt{P}$. The points obtained in this way have the
following properties: the vector $\Phi_s$ is of the form 
$\Phi^{\otimes r/2}$, $\Phi$ factors as 
$$
\xymatrix{
{V^{\otimes 2}\otimes \SO_X(-2m)} \ar[rr]^{\Phi} \ar[d] & & 
{H^0(\SO_X(2m))}\otimes \SO_X(-2m) \ar[d] \\
{E^{\otimes 2}} \ar[rr]^{\varphi} & & {\SO_X \;,}
}
$$
the homomorphism $\Psi$ factors as
$$
\xymatrix{
{\bigwedge{}^r V \otimes \SO_X(-rm)} \ar[rr]^{\Psi} \ar[d] & & 
{H^0(\SO_X(rm))\otimes \SO_X(-rm)} \ar[d] \\
{\det E} \ar[rr]^{\psi} & & {\SO_X \;,}
}
$$
and $\det(\phi)=\psi^2$ as in (SOS5).

Let $\wt{Z}'$ be the closed subset of $\SH\times \wt{P}$
defined by these properties. Given a point $z\in \wt{Z}'$
we can recover the tuple up to isomorphism.
Define the parameter space $\wt{Z}$ as the closure in $\wt{Z}'$ 
of those points obtained
from semistable special orthogonal sheaves.

Let $\delta$ be a polynomial as in (\ref{delta}) and with $\delta_1>0$. 
Define a polarization on $\wt{Z}$ by 
$$
\SO_{\wt{Z}}(n_1,n_2):= p^*_\SH\SO^{}_\SH(n_1) \otimes 
p^*_P \SO^{}_P(\frac{2 n_2}{r}),
$$
where $n_2$ is a multiple of $r/2$, $n_1$ is an integer, and
$$
\frac{n_2}{n_1}=\frac{P(l)\delta(m)-\delta(l)P(m)}{P(m)-2\delta(m)}
$$
The projective scheme $\wt{Z}$ is preserved by the natural $\slv$ action,
and this action has a natural linearization on
$\SO_{\wt{Z}}(n_1,n_2)$.

\begin{proposition}
A point $(g,E,\varphi,\psi)$ is GIT-(semi)stable if and only if
the special orthogonal sheaf $(E,\varphi,\psi)$ is (semi)stable
(definition \ref{stsor}).
\end{proposition}

\begin{proof}
The parameter space $Z$ for orthogonal sheaves is a subscheme
of $\SH \times P$, where 
$$
P=\PP\big((V^{\otimes 2})^\vee \otimes H^0(\SO_X(2m))\big)
$$
(this is a particular case of the parameter space defined in 
section \ref{secgitconstruction}). Let $\SO_Z(n_1,n_2)$ be 
the polarization defined in (\ref{son1n2}), and consider the
natural linearization of the action of $\slv$ on this polarization.
There is a morphism 
$$
\begin{array}{ccc}
{\wt{Z}} & \stackrel{f}{\too} & {Z} \\ 
{(g,E,\varphi,\psi)} & \longmapsto & {(g,E,\varphi)}
\end{array}
$$
with $f^* \SO_Z(n_1,n_2)=\SO_{\wt Z}(n_1,n_2)$. This morphism is
equivariant with respect to $\slv$, and the linearizations are compatible.
Property (SOS5$'$) implies that $f$ is finite \'etale 
(because $\wt{Z}$ is given
locally by the equation $\det(\varphi(x))=\psi(x)^2$),  
and then it follows that 
a point in $\wt{Z}$ is GIT-(semi)stable
if and only if its image in $Z$ is GIT-(semi)stable. The result
follows from theorem \ref{thm3.2},
proposition \ref{critGM}, and definition \ref{stsor}.
\end{proof}

Let $\SM_{\sor}$ be the functor of families of 
semistable special orthogonal sheaves.
Let $\FM_{\sor}$ be the GIT quotient of $\wt Z$ by $\slv$.
Let $(E,\varphi,\psi)$ be a semistable special orthogonal sheaf. 
Let $E^S$ and
$\varphi^S$ be defined as in section \ref{sectheorem}. There is a
natural isomorphism between $\det E^S $ and $\det E $, then composing
with $\psi$ we obtain an isomorphism $\psi^S:\det E^S \to \SO_X$.

Let $(E,\varphi,\psi)$ and $(E',\varphi',\psi')$ be two 
semistable special orthogonal sheaves. They are $S$-equivalent 
if and only if $(E^S,\varphi^S,\psi^S)$ is isomorphic to 
$({E'}^S,{\varphi'}^S, {\psi'}^S)$.

\begin{theorem}
\label{modsor}
The projective scheme ${\FM}_{\sor}$ 
is the coarse moduli space of S-equivalence classes of 
semistable special orthogonal sheaves.
There is an open subset $\FM_{\sor}^0$ corresponding to semistable
special orthogonal bundles. 
\end{theorem}

The proof is completely analogous to the proof of theorem
\ref{mainthm} (section \ref{sectheorem}).

\section{$\glr$-representational pairs}
\label{secrelatedmodulis}

Once we have constructed the moduli space of \tensors, it is easy
to obtain moduli spaces for $\glr$-representational pairs. 
In the case of
$\dim(X)=1$, this is done in \cite{Sch}, but since it does not 
depend on the dimension of the base $X$, the same arguments apply   
here. 
In \cite{Ba}, Banfield considered pairs $(P,\sigma)$, where
$P$ is a principal $G$-bundle ($G$ any reductive group),
and $\sigma$ is a section
associated to $P$ by a fixed representation $\rho$.
He defined stability, and proved a Hitchin-Kobayashi correspondence.
Now we will construct the moduli space, when $G=\glr$.

Fix a polynomial $\delta$ as in (\ref{delta}).
Let $\rho:\glr \to \gln$ be a representation sending the 
center of $\glr$ to the center of $\gln$. Consider a triple
\begin{equation}
\label{rep}
(E,\; \psi:E_\rho\to \fixed_u,\; u),
\end{equation}
where $E$ is a vector
bundle of rank $r$ on $X$, and $E_\rho$ is the vector bundle of 
rank $n$ associated to $E$ and $\rho$. Using \cite[prop. 15.47]{F-H},
it can be shown that there exist integers $s>0$, $b$, $c>0$, 
and a vector bundle $F$ such
that 
$$
(E^{\otimes s})^{\oplus c}\otimes (\det E)^{-\otimes b}
\; \cong \; E_\rho \oplus F
$$
(see \cite[cor 1.1.2]{Sch} for details). 
Then a triple 
(\ref{rep}) is equivalent to a \tensor\ $(E,\varphi,u)$
such that 
\begin{equation}
\label{f0}
\varphi|^{}_F=0,
\end{equation}
and we say that the triple is
$\delta$-(semi)stable if the corresponding \tensor\ is. 
Since the condition (\ref{f0}) is closed, the moduli space of 
$\delta$-semistable triples is a closed subscheme 
of $\FM^0_\delta$,
the open subscheme corresponding to \tensors\ with
$E$ locally free.
It is easy to check that the definition of stability in \cite{Ba}
coincides with our slope-$\tau$-stability.

In \cite{MR}, Mundet generalized Banfield's work.
He fixes a Kaehler manifold $Y$ and an action $\rho$ of a reductive
group $G$ on $Y$, and considers
pairs $(P,\sigma)$,  where $\sigma$ is a section of the 
associated fiber bundle $P\times_G Y$. He defined stability and
proved a Hitchin-Kobayashi correspondence. Now we will construct
the moduli space, for the case when $G=\glr$, and $Y$ is a 
projective (or more generally, quasi-projective) scheme. 

Consider an action of $\glr$ on a projective scheme $Y$, together
with a linearization of the action on an ample line bundle $L$ on $Y$.
Assume that the center of $\glr$ acts trivially on $Y$.
Consider a pair 
\begin{equation}
\label{proj}
(P,\; \sigma:X\too P\times_{\glr} Y),
\end{equation}
where $P$ is a principal $\glr$ bundle on $X$, and $\sigma$ is
a section of the fiber bundle associated to $P$ with fiber $Y$.
We fix the topology type of $P$ and the homology class 
$[\sigma(X)]$ of the image.
Fix $k$ large enough so that we have a natural embedding
$F\inj \PP(H^0(F,L^{\otimes k})^\vee)$.
Since the action of the center of $\glr$ is trivial on $Y$,
the induced representation 
$$
\rho:\glr \too \operatorname{GL}(H^0(F,L^{\otimes k}))
$$
sends the center of $\glr$ to the center of 
$\operatorname{GL}(H^0(F,L^{\otimes k}))$.
Let $E$ be the rank $r$ vector bundle corresponding to $P$.
Since we have fixed the topology type of $P$, the Hilbert
polynomial $P_E$ of $E$ is also fixed.
The section $\sigma$ gives a homomorphism
$$
\psi: E_\rho \too \fixed_u,
$$
for some line bundle $\fixed_u$, whose degree $a$ depends on the
homology class $[\sigma(X)]$ 
of the image. Take $\param=\Pic^a(X)$, and let $\sfixed$ be a 
Poincare bundle.
We obtain that a pair (\ref{proj}) is equivalent
to a triple (\ref{rep}) with the property that
the section $\psi':X\to \PP(E_\rho^\vee)$
factors through $P\times_{\glr} Y$. We define a pair
(\ref{proj}) to be $\delta$-semistable if the corresponding
triple is, and hence the moduli space of $\delta$-semistable 
pairs (\ref{proj}) is a closed subscheme of $\FM^0_\delta$. 
We can also take $Y$ to be
quasi-projective, and the moduli space will also be a subscheme
(not necessarily closed) of $\FM^0_\delta$.

\end{document}